# LARGE DEVIATIONS FOR TRAPPED INTERACTING BROWNIAN PARTICLES AND PATHS[1]

BY STEFAN ADAMS, JEAN-BERNARD BRU AND WOLFGANG KÖNIG

*Max-Planck Institute for Mathematics in the Sciences and Dublin Institute for Advanced Studies, Johannes-Gutenberg-Universität Mainz and Universität Leipzig*

We introduce two probabilistic models for $N$ interacting Brownian motions moving in a trap in $\mathbb{R}^d$ under mutually repellent forces. The two models are defined in terms of transformed path measures on finite time intervals under a trap Hamiltonian and two respective pair-interaction Hamiltonians. The first pair interaction exhibits a *particle* repellency, while the second one imposes a *path* repellency.

We analyze both models in the limit of diverging time with fixed number $N$ of Brownian motions. In particular, we prove large deviations principles for the normalized occupation measures. The minimizers of the rate functions are related to a certain associated operator, the Hamilton operator for a system of $N$ interacting trapped particles. More precisely, in the particle-repellency model, the minimizer is its ground state, and in the path-repellency model, the minimizers are its ground product-states. In the case of path-repellency, we also discuss the case of a Dirac-type interaction, which is rigorously defined in terms of Brownian intersection local times. We prove a large-deviation result for a discrete variant of the model.

This study is a contribution to the search for a mathematical formulation of the quantum system of $N$ trapped interacting bosons as a model for *Bose–Einstein condensation*, motivated by the success of the famous 1995 experiments. Recently, Lieb et al. described the large-$N$ behavior of the ground state in terms of the well-known *Gross–Pitaevskii* formula, involving the scattering length of the pair potential. We prove that the large-$N$ behavior of the ground product-states is also described by the Gross–Pitaevskii formula, however, with the scattering length of the pair potential replaced by its integral.

Received December 2004; revised January 2006.

[1]Supported in part by DFG Priority Research Program DE 663/1-3 and DFG Grant AD 194/1-1.

## 1. Introduction and results.

1.1. *Introduction and motivation.*   Consider a large number of identical quantum particles in a trap under the influence of a mutually repellent pair interaction. We assume that their wave function is invariant under permutations of the single-particle variables. In physics particle ensembles whose many-body wave functions have this invariance property are called *boson systems.* A specific property of bosons is that, for large number of particles, at extremely low temperature, they undergo a phase transition, the so called *Bose–Einstein condensation.* This formally means that a macroscopic portion, the condensate, is described by one *single-particle* wave function. Bose–Einstein condensation was theoretically predicted by S. N. Bose and A. Einstein in the 1920s, and a huge theoretical literature has been accumulated since. After appropriate cooling methods had been developed, the first experimental realization of this condensation succeeded in 1995 [4, 6, 10]. For this remarkable achievement, the Nobel prize in physics 2001 was awarded to E. A. Cornell, W. Ketterle and C. E. Wieman. A comprehensive account on Bose–Einstein condensation is the recent monograph [28].

Motivated by the experimental success, in a series of papers Lieb, Seiringer and Yngvason [22, 23, 24, 25] obtained a mathematical foundation of Bose–Einstein condensation at zero temperature. The mathematical formulation of the $N$-particle boson system is in terms of an $N$-particle Hamilton operator, $\mathcal{H}_N$, whose ground states describe the bosons under the influence of a trap potential and a pair potential, see (1.7). Lieb et al. rigorously proved that the ground state energy per particle of $\mathcal{H}_N$ (after proper rescaling of the pair potential) converges toward the energy of the well-known *Gross–Pitaevskii* functional, and the ground state is approximated by the $N$-fold product of the Gross–Pitaevskii minimizer. Moreover, they also showed the convergence of the reduced density matrix, which implies the Bose–Einstein condensation. As had been generally predicted, the *scattering length* of the pair interaction potential plays a key role in this description.

Much thermodynamic information about the boson system is contained in the traces of the Boltzmann factor $e^{-\beta \mathcal{H}_N}$ for $\beta > 0$, like the free energy, or the pressure. Since the 1960s, interacting Brownian motions are generally used for probabilistic representations for these traces. The parameter $\beta$, which is interpreted as the inverse temperature of the system, is then the length of the time interval of the Brownian motions. However, the traces do not contain much information about the ground *state.* Since the pioneering work of Donsker and Varadhan in the early 1970s, it is basically known that the ground states are intimately linked with the Brownian *occupation measures.* This link is established via the theory of *large deviations* for diverging time, which corresponds to vanishing temperature.

In the present article we study models of a fixed number $N$ of trapped interacting Brownian motions given in terms of transformed measures for paths of length $\beta$. A *trap Hamiltonian*, putting "hard" or "soft" walls, keeps the motions in a bounded region. We consider basically two different types of *pair-interaction Hamiltonians*, both imposing mutually repellent interaction: the first one, which we call the *canonical ensemble* model, imposes *particle*-repellency, while the second, which we call the *Hartree* model, imposes *path*-repellency. More precisely, in the canonical ensemble model the $N$ motions interact with each other at common time units, and in the Hartree model, the paths interact with the mean of the paths of the other motions. The canonical ensemble model comprises the above mentioned standard representation of the traces of $e^{-\beta \mathcal{H}_N}$. The Hartree model is apparently introduced here for the first time; it cannot be represented in terms of the trace of any Boltzmann factor.

It is the main aim of this paper to study the large-$\beta$ behavior of the canonical ensemble and the Hartree model in terms of large deviations principles for the Brownian occupation measures. One of our main results is a large deviation principle for the joint occupation measure of the $N$-tuple of the motions in the canonical ensemble model (Theorem 1.5). The rate function is given by the energy of the operator $\mathcal{H}_N$. In particular, we prove a law of large numbers for the occupation measure toward the ground state. The large-$\beta$ behavior of the Hartree model is different and is described in terms of the ground *product states* of $\mathcal{H}_N$, that is, the minimizers of the energy of $\mathcal{H}_N$ among all product states. Our main result here is a large deviation principle for the $N$-tuple of the occupation measures of the motions (Theorem 1.7). The rate function is the $\mathcal{H}_N$-energy of the product of the components of the tuple. In particular, we have a law of large numbers for the occupation measure tuple toward the set of the ground product states. We also discuss a mathematical idealization of the Hartree model, where the pair-interaction potential is replaced by the Dirac-measure at zero. This model (which is trivial for the canonical ensemble model) is defined in terms of Brownian *intersection local times*, an object whose large deviations properties are currently much studied from a probabilistic point of view. We prove analogous large deviation results for a discrete variant of the model (Theorem 1.12). Finally, we also study the ground product states of $\mathcal{H}_N$ in the limit $N \to \infty$ (Theorem 1.14). In dimension $d = 3$, the rescaling of the pair potential is the same as was considered by Lieb et al. It turns out that the energy converges toward the Gross–Pitaevskii formula, and the mean converges toward the minimizer. Instead of the scattering length of the pair-interaction potential, its integral, which is a strictly larger number, is the decisive parameter. In dimension $d = 2$, we prove the same result, however, the scaling is here different from the one considered by Lieb et al.

The replacement of the ground state energy by the ground product state energy is known as the *Hartree–Fock approach* [12]. One of the main novelties of the present article is a probabilistic version of the Hartree–Fock idea for interacting Brownian motions; indeed, the product states are interpreted in terms of the occupation measures in the Hartree model. Another main novelty is a rigorous proof of the asymptotic relation, for diverging particle number, between the ground product states and the Gross–Pitaevskii minimizer with parameter equal to the integral of the pair-interaction potential. This relation and the one between (unrestricted) ground states and the scattering length is a general issue and has been phenomenologically discussed elsewhere [14, 29]. However, the former has previously not been rigorously proved for the $N$-body problem.

The present paper provides the mathematical basis for the rigorous probabilistic analysis of interacting Brownian motions as models for large quantum particle systems at positive temperature. In [1] we give an asymptotic analysis of the Hartree model at fixed finite time for diverging particle number. Future work will be devoted to the same question for the canonical ensemble (i.e., for the trace of $e^{-\beta \mathcal{H}_N}$), which is important from the physical point of view, but also more ambitious. Specific permutation symmetric systems of Brownian bridges, characteristic for boson systems, will be analyzed in [2]. The latter hopefully leads to a probabilistic model for Bose–Einstein condensation at *fixed positive* temperature.

The remainder of this section is structured as follows. We first define the two basic variational problems in Section 1.2. In Section 1.3 we define the canonical ensemble model and the Hartree model. Our results on the large deviations behaviors of these two models are also presented in Section 1.3. In Section 1.4 we discuss the Dirac-type interaction. The large-$N$ behavior of our key variational formulas is treated in Section 1.5. The notion of a large-deviation principle and of the scattering length are recalled in Section 1.6.

The subsequent sections are devoted to the proofs: in Sections 2, 3 and 5 we prove our results on the Hartree model, its Dirac-type variant and the canonical ensemble model, respectively, and in Section 4 we analyze the ground product states and its large-$N$ behavior. Throughout the paper we consider only dimensions $d \geq 2$; in Section 1.5 we restrict to the case $d \in \{2, 3\}$.

1.2. *The variational problems.* We introduce two fundamental variational problems for the energy of $N$ particles in $\mathbb{R}^d$ under the influence of two potentials.

1.2.1. *The potentials.* Our two fundamental ingredients are a *trap potential*, $W$, and a *pair-interaction potential*, $v$. Our assumptions on $W$ are the

following:

(1.1)
$$W : \mathbb{R}^d \to [0, \infty] \text{ is measurable and locally integrable on } \{W < \infty\}$$
$$\text{with } \liminf_{R \to \infty} \inf_{|x| > R} W(x) = \infty.$$

In order to avoid trivialities, we assume that $\{W < \infty\}$ is either equal to $\mathbb{R}^d$ or is a bounded connected open set containing the origin.

Our assumptions on $v$ are the following. By $B_r(x)$ we denote the open ball with radius $r$ around $x \in \mathbb{R}^d$:

(1.2)
$$v : [0, \infty) \to \mathbb{R} \cup \{+\infty\} \text{ is measurable and bounded from below,}$$
$$a := \sup\{r \geq 0 : v(r) = \infty\} \in [0, \infty), \qquad v|_{[\eta, \infty)} \text{ is bounded } \forall \eta > a.$$

Note that we also admit $v(a) = +\infty$. We are mainly interested in the case where $v$ has a singularity, that is, either $a > 0$, or $a = 0$ and $\lim_{r \downarrow 0} v(r) = \infty$. Examples include also super-stable potentials and potentials of Lennard–Jones type [30]. According to integrability properties near the origin, we distinguish two different classes as follows.

DEFINITION 1.1. We call the interaction potential $v$ a *soft-core* potential if $a = 0$ and $\int_{B_1(0)} v(|x|) \, dx < +\infty$. Otherwise [i.e., if $a > 0$, or if $a = 0$ and $\int_{B_1(0)} v(|x|) \, dx = +\infty$], we call the interaction potential a *hard-core* potential.

We shall need the following $dN$-dimensional versions of the trap and the interaction potential:

(1.3)
$$\mathfrak{W}(x) = \sum_{i=1}^{N} W(x_i) \quad \text{and} \quad \mathfrak{v}(x) = \sum_{1 \leq i < j \leq N} v(|x_i - x_j|),$$

where $x = (x_1, \ldots, x_N) \in \mathbb{R}^{dN}$. Our potential $\mathfrak{W} + \mathfrak{v}$ is locally integrable precisely on the set

(1.4)
$$\Omega = \{W < \infty\}^N \cap \begin{cases} \mathbb{R}^{dN}, & \text{if } v \text{ is soft-core,} \\ U_a, & \text{if } v \text{ is hard-core,} \end{cases}$$

where

(1.5)
$$U_\eta := \bigcap_{1 \leq i < j \leq N} \{x \in \mathbb{R}^{dN} : |x_i - x_j| > \eta\}, \qquad \eta \geq 0.$$

Note that $U_\eta$ is a connected set in $d \geq 2$.

1.2.2. *The ground state energy.* Our first fundamental object is the ground-state energy per particle

$$
\begin{aligned}
(1.6) \quad \chi_N &= \frac{1}{N} \inf_{h \in H^1(\Omega) \,:\, \|h\|_2 = 1} \langle h, \mathcal{H}_N h \rangle \\
&= \frac{1}{N} \inf_{h \in H^1(\Omega) \,:\, \|h\|_2 = 1} [\|\nabla h\|_2^2 + \langle \mathfrak{W}, h^2 \rangle + \langle \mathfrak{v}, h^2 \rangle],
\end{aligned}
$$

of the $N$-particle Hamilton operator

$$
(1.7) \qquad \mathcal{H}_N = -\Delta + \mathfrak{W} + \mathfrak{v} \qquad \text{on } L^2(\Omega).
$$

Here $H^1(\Omega) = \{f \in L^2(\Omega) : \nabla f \in L^2(\Omega)\}$ is the usual Sobolev space, and $\nabla$ is the distributional gradient. We summarize some facts about the ground states of $\mathcal{H}_N$.

LEMMA 1.2 (Ground states of $\mathcal{H}_N$). *Fix $N \in \mathbb{N}$.*

(i) *There is a unique minimizer $h_* \in H^1(\Omega)$ on the right-hand side of* (1.6).

(ii) *The minimizer $h_*$ satisfies the variational equation*

$$
(1.8) \qquad \Delta h_* = \mathfrak{W} h_* + \mathfrak{v} h_* - N \chi_N h_*.
$$

(iii) *The minimizer $h_*$ is positive everywhere on $\Omega$ and continuously differentiable, and its first partial derivatives are $\alpha$-Hölder continuous for any $\alpha < 1$.*

PROOF. The proof is standard and is therefore omitted. See the proof of Lemma 1.3 below or Appendix A of [22] for the treatment of similar problems. □

1.2.3. *The ground product state energy.* Introduce the ground product state energy of $\mathcal{H}_N$, that is,

$$
(1.9) \quad \chi_N^{(\otimes)} = \frac{1}{N} \inf_{h_1, \ldots, h_N \in H^1(\mathbb{R}^d) \,:\, \|h_i\|_2 = 1 \,\forall i} \langle h_1 \otimes \cdots \otimes h_N, \mathcal{H}_N h_1 \otimes \cdots \otimes h_N \rangle.
$$

The replacement of the ground state energy, $\chi_N$, by the ground product state energy, $\chi_N^{(\otimes)}$, is known as the *Hartree–Fock approach* [12]. Sometimes, the formula in (1.9) is called the *Hartree formula*. Obviously,

$$
(1.10) \qquad \chi_N^{(\otimes)} \geq \chi_N.
$$

We can also write

$$
(1.11) \quad \begin{aligned}
\chi_N^{(\otimes)} = \frac{1}{N} \inf_{h_1, \ldots, h_N \in H^1(\mathbb{R}^d) \,:\, \|h_i\|_2 = 1 \,\forall i} \Bigg[ &\sum_{i=1}^N [\|\nabla h_i\|_2^2 + \langle W, h_i^2 \rangle] \\
&+ \sum_{1 \leq i < j \leq N} \langle h_i^2, V h_j^2 \rangle \Bigg],
\end{aligned}
$$

where $V$ denotes the integral operator with kernel $v \circ |\cdot|$, either defined for functions by $Vf(x) = \int_{\mathbb{R}^d} v(|x - y|) f(y)\, dy$ or for measures by $V\mu(x) = \int_{\mathbb{R}^d} \mu(dy) v(|x - y|)$.

The main assertions on the formula in (1.9) and its minimizers are summarized as follows.

LEMMA 1.3 (Product ground states of $\mathcal{H}_N$). *Fix $N \in \mathbb{N}$.*

(i) *There exists at least one minimizer $(h_1, \ldots, h_N)$ of the right-hand side of* (1.11). *The set of minimizers is compact and invariant under permutation of the functions $h_1, \ldots, h_N$.*

(ii) *Any minimizer $(h_1, \ldots, h_N)$ satisfies the system of differential equations*

$$(1.12) \qquad \Delta h_i = -\lambda_i h_i + W h_i + h_i \sum_{j \neq i} V h_j^2, \qquad i = 1, \ldots, N,$$

*with $\lambda_i = \|\nabla h_i\|_2^2 + \langle W, h_i^2 \rangle + \sum_{j \neq i} \langle h_i^2, V h_j^2 \rangle$. Furthermore, $\|h_i\|_\infty \leq C_d(\lambda_i - (N-1)\inf v)^{d/4}$ for any $i \in \{1, \ldots, N\}$, where $C_d > 0$ depends on the dimension $d$ only.*

(iii) *Let $v$ be soft-core, assume that $d \in \{2, 3\}$, and let $(h_1, \ldots, h_N)$ be any minimizer. Assume that $v|_{(0,\eta)} \geq 0$ for some $\eta > 0$. In $d = 3$, furthermore assume that*

$$(1.13) \qquad \int_{B_1(0)} |v(|y|)|^{1+\delta}\, dy < \infty \qquad \text{for some } \delta > 0.$$

*Then every $h_i$ is positive everywhere in $\mathbb{R}^d$ and continuously differentiable, and all first partial derivatives are $\alpha$-Hölder continuous for any $\alpha < 1$.*

(iv) *Let $v$ be hard-core, assume that $d \in \{2, 3\}$, and let $(h_1, \ldots, h_N)$ be any minimizer. Then every $h_i$ is continuously differentiable in the interior of its support, and all first partial derivatives are $\alpha$-Hölder continuous for any $\alpha < 1$.*

For the proof see Sections 4.1–4.4.

REMARK 1.4. (i) Unlike for the ground states of $\mathcal{H}_N$ in (1.6), there is no convexity argument available for the formula in (1.9). This is due to the fact that a convex combination of tensor-products of functions is not tensor-product in general, and hence, the domain of the infimum in (1.9) is not a convex subset of $H^1(\mathbb{R}^{dN})$. However, for $h_2, \ldots, h_N$ fixed, the minimization over $h_1$ enjoys the analogous convexity properties on $H^1(\mathbb{R}^d)$ as the minimization in (1.6).

(ii) If $v$ is hard-core, it is easy to see that the distances between the supports of $h_1, \ldots, h_N$ have to be no smaller than $a$ [see (1.2)] in order to make the value of $\langle h_1 \otimes \cdots \otimes h_N, \mathcal{H}_N h_1 \otimes \cdots \otimes h_N \rangle$ finite. The potential $\sum_{j \neq i} V h_j^2$ is equal to $\infty$ in the $a$-neighborhood of the union of the supports of $h_j$ with $j \neq i$, and $h_i$ is equal to zero there (we regard $0 \cdot \infty$ as 0). In particular, minimizers of (1.9) are not of the form $(h, \ldots, h)$. This shows that the inequality in (1.10) is even strict. In the soft-core case, this statement is not obvious at all. A partial result on this question in $d = 3$ will be a by-product of Section 1.5 below.

1.3. *The Brownian models.* We introduce two different models of interacting Brownian motions. These models are given in terms of transformed measures for paths of length $\beta$ in terms of certain Hamiltonians. Let a family of $N$ independent Brownian motions, $(B_t^{(1)})_{t \geq 0}, \ldots, (B_t^{(N)})_{t \geq 0}$, in $\mathbb{R}^d$ with generator $-\Delta$ be given. The Hamiltonians of both models possess a trap part and a pair-interaction part. The trap part is for both models the same, namely,

$$(1.14) \qquad H_{N,\beta} = \sum_{i=1}^{N} \int_0^\beta W(B_s^{(i)}) \, ds.$$

We assume that the joint starting distribution of the motions, that is, the distribution of the vector $(B_0^{(1)}, \ldots, B_0^{(N)})$, is concentrated on a compact subset of $\Omega$.

1.3.1. *The canonical ensemble model.* The Hamiltonian of our first model consists of two parts: the trap part given in (1.14), and a pair-interaction part,

$$(1.15) \qquad G_{N,\beta} = \sum_{1 \leq i < j \leq N} \int_0^\beta v(|B_s^{(i)} - B_s^{(j)}|) \, ds.$$

We look at the distribution of the $N$ Brownian motions under the transformed path measure

$$
(1.16) \qquad d\widehat{\mathbb{P}}_{N,\beta} = \frac{1}{Z_{N,\beta}} \exp(-H_{N,\beta} - G_{N,\beta}) \, d\mathbb{P},
$$
$$
\text{where } Z_{N,\beta} = \mathbb{E}(\exp(-H_{N,\beta} - G_{N,\beta})).
$$

We call $\widehat{\mathbb{P}}_{N,\beta}$ the *canonical ensemble model*, since it is derived, via a Feynman–Kac formula, from the trace-class operator of the canonical ensemble, $e^{-\beta \mathcal{H}_N}$ (see Remark 1.6 below). It is a model for $N$ Brownian motions in a trap $W$ with the presence of a repellent pair interaction. We can conceive the $N$-tuple of the motions, $B_t = (B_t^{(1)}, \ldots, B_t^{(N)})$, as one Brownian motion in

$\mathbb{R}^{dN}$. Introduce the normalized *occupation measure* of the $dN$-dimensional motion,

$$(1.17) \qquad \mu_\beta(dx) = \frac{1}{\beta} \int_0^\beta \delta_{B_s}(dx)\, ds,$$

which is a random element of the set $\mathcal{M}_1(\mathbb{R}^{dN})$ of probability measures on $\mathbb{R}^{dN}$. It measures the time spent by the tuple of $N$ Brownian motions in a given region. Note that there is only one time scale involved for all the motions, that is, the Brownian particles interact with each other at common time units. We can write the Hamiltonians in terms of the occupation measure as

$$(1.18) \qquad H_{N,\beta} = \beta\langle \mathfrak{W}, \mu_\beta\rangle \quad \text{and} \quad G_{N,\beta} = \beta\langle \mathfrak{v}, \mu_\beta\rangle,$$

where the functions $\mathfrak{W}, \mathfrak{v} \colon \mathbb{R}^{dN} \to \mathbb{R}$ are introduced in (1.3).

It turns out that the large-$\beta$ behavior of the canonical ensemble model is described by the ground state of the operator $\mathcal{H}_N$ in (1.7). In Section 1.6.2 below we recall the notion of a principle of large deviations. The rate function $I_N$ appearing in Theorem 1.5 is the well-known *Donsker–Varadhan rate function* on $\mathbb{R}^{dN}$ defined by

$$(1.19) \qquad I_N(\mu) = \begin{cases} \left\| \nabla \sqrt{\dfrac{d\mu}{dx}} \right\|_2^2, & \text{if } \sqrt{\dfrac{d\mu}{dx}} \in H^1(\mathbb{R}^{dN}) \text{ exists,} \\ \infty, & \text{otherwise.} \end{cases}$$

Note that the energy functional $\langle h, \mathcal{H}_N h\rangle$ may be rewritten $\langle h, \mathcal{H}_N h\rangle = I_N(\mu) + \langle \mathfrak{W}, \mu\rangle + \langle \mathfrak{v}, \mu\rangle$ for the probability measure $\mu(dx) = h^2(x)\, dx$.

THEOREM 1.5 (Canonical ensemble model at late times). *Fix $N \in \mathbb{N}$.*

(i)

$$(1.20) \qquad \lim_{\beta \to \infty} \frac{1}{N\beta} \log \mathbb{E}(\exp(-H_{N,\beta} - G_{N,\beta})) = -\chi_N,$$

*where $\chi_N$ is the ground-state energy per particle of the $N$-particle operator $\mathcal{H}_N$ given in (1.6).*

(ii) *As $\beta \to \infty$, the distribution of $\mu_\beta$ on $\mathcal{M}_1(\mathbb{R}^{dN})$ under $\widehat{\mathbb{P}}_{N,\beta}$ satisfies a principle of large deviation with speed $\beta$ and rate function $I_N$ given by*

$$(1.21) \quad I_N(\mu) = I_N(\mu) + \langle \mathfrak{W}, \mu\rangle + \langle \mathfrak{v}, \mu\rangle - N\chi_N \qquad \text{for } \mu \in \mathcal{M}_1(\mathbb{R}^{dN}).$$

(iii) *The distribution of $\mu_\beta$ under $\widehat{\mathbb{P}}_{N,\beta}$ converges weakly toward the measure $h_*(x)^2\, dx$, where $h_*$ is the unique minimizer in (1.6).*

For the proof see Section 5.

REMARK 1.6. It is well known [18] that the bottom of the spectrum of $\mathcal{H}_N$ is related to the large-$\beta$ behavior of the trace of $e^{-\beta\mathcal{H}_N}$, more precisely,

$$(1.22) \qquad \chi_N = -\lim_{\beta\to\infty}\frac{1}{N\beta}\log\mathrm{Tr}(e^{-\beta\mathcal{H}_N}).$$

Using the Feynman–Kac formula for traces [7, 18], we have

$$(1.23) \qquad \mathrm{Tr}(e^{-\beta\mathcal{H}_N}) = \int_\Omega dx\, \mathbb{E}_x^{(\beta)}(e^{-\beta\langle\mathfrak{W}+\mathfrak{v},\mu_\beta\rangle}),$$

where $\mathbb{E}_x^{(\beta)}$ is the expectation with respect to a Brownian bridge of length $\beta$ that starts and ends at $x\in\mathbb{R}^{dN}$. Note the close relation to Theorem 1.5.

1.3.2. *The Hartree model.* Our second Brownian model is defined in terms of another Hamiltonian. We keep the trap Hamiltonian $H_{N,\beta}$ as in (1.14), but the interaction Hamiltonian is now

$$(1.24) \qquad K_{N,\beta} = \sum_{1\le i<j\le N}\frac{1}{\beta}\int_0^\beta\int_0^\beta v(|B_s^{(i)}-B_t^{(j)}|)\,ds\,dt.$$

Note that the $i$th Brownian motion interacts with the mean of the whole path of the $j$th motion, taken over all times before $\beta$. Hence, the interaction is not a *particle* interaction, but a *path* interaction. The interaction (1.24) is related to Polaron-type models [5, 13], where instead of several paths a single path is considered. We consider the corresponding transformed path measure,

$$(1.25) \qquad d\widehat{\mathbb{P}}_{N,\beta}^{(\otimes)} = \frac{1}{Z_{N,\beta}^{(\otimes)}}\exp(-H_{N,\beta}-K_{N,\beta})\,d\mathbb{P},$$
$$\text{where } Z_{N,\beta}^{(\otimes)} = \mathbb{E}(\exp(-H_{N,\beta}-K_{N,\beta})).$$

In Theorem 1.7 below it turns out that the large-$\beta$ behavior of $Z_{N,\beta}^{(\otimes)}$ is intimately related to the Hartree formula in (1.9). Therefore, we call this model the *Hartree model*. At the end of this section we comment on its physical relevance.

We introduce the normalized occupation measure of the $i$th motion,

$$(1.26) \qquad \mu_\beta^{(i)}(dx) = \frac{1}{\beta}\int_0^\beta \delta_{B_s^{(i)}}(dx)\,ds \in \mathcal{M}_1(\mathbb{R}^d).$$

The tuple of the $N$ occupation measures, $(\mu_\beta^{(1)},\ldots,\mu_\beta^{(N)})$, plays a particular role in this model. We can write the Hamiltonians as

$$(1.27) \qquad \begin{aligned} H_{N,\beta} &= \beta\sum_{i=1}^N\langle W,\mu_\beta^{(i)}\rangle = \beta\langle\mathfrak{W},\mu_\beta^{\otimes}\rangle \quad\text{and}\\ K_{N,\beta} &= \beta\sum_{1\le i<j\le N}\langle\mu_\beta^{(i)},V\mu_\beta^{(j)}\rangle = \beta\langle\mathfrak{v},\mu_\beta^{\otimes}\rangle, \end{aligned}$$

where we recall (1.3) and the operator $V$ with kernel $v \circ |\cdot|$, and $\mu_\beta^\otimes = \mu_\beta^{(1)} \otimes \cdots \otimes \mu_\beta^{(N)}$ is the product measure. Now we formulate the main result on the large-time behavior of the Hartree model.

THEOREM 1.7 (Hartree model at late times). *Assume that $W$ and $v$ are continuous in $\{W < \infty\}$, respectively, in $\{v < \infty\}$. Furthermore, assume in the soft-core case that there exists an $\varepsilon > 0$ and a decreasing function $\tilde{v} : (0, \varepsilon) \to \mathbb{R}$ with $v \leq \tilde{v}$ on $(0, \varepsilon)$, which satisfies $\int_{B_\varepsilon(0)} G(0, y)\tilde{v}(|y|)\, dy < \infty$, where $G$ denotes the Green's function of the free Brownian motion on $\mathbb{R}^d$. Fix $N \in \mathbb{N}$.*

(i) *For $\chi_N^{(\otimes)}$ as defined in (1.9),*

$$\lim_{\beta \to \infty} \frac{1}{N\beta} \log \mathbb{E}(\exp(-H_{N,\beta} - K_{N,\beta})) = -\chi_N^{(\otimes)}. \tag{1.28}$$

(ii) *As $\beta \to \infty$, the distribution of the tuple $(\mu_\beta^{(1)}, \ldots, \mu_\beta^{(N)})$ of Brownian occupation measures on $\mathcal{M}_1(\mathbb{R}^d)^N$ under $\widehat{\mathbb{P}}_{N,\beta}^{(\otimes)}$ satisfies a large deviation principle with speed $\beta$ and rate function*

$$I_N^{(\otimes)}(\mu_1, \ldots, \mu_N) = \sum_{i=1}^N I_1(\mu_i) + \langle \mathfrak{W}, \mu^\otimes \rangle + \langle \mathfrak{v}, \mu^\otimes \rangle - N\chi_N^{(\otimes)},$$
$$\mu_1, \ldots, \mu_N \in \mathcal{M}_1(\mathbb{R}^d), \tag{1.29}$$

*where $I_1$ is defined in (1.19), and $\mu^\otimes = \mu_1 \otimes \cdots \otimes \mu_N$ is the product measure.*

(iii) *The distribution of $(\mu_\beta^{(1)}, \ldots, \mu_\beta^{(N)})$ under $\widehat{\mathbb{P}}_{N,\beta}^{(\otimes)}$ is attracted by the set of minimizers in (1.11).*

For the proof see Section 2.

REMARK 1.8. (1) The additional assumption of continuity of $W$ and $v$ is necessary only in the proof of the upper bound, where we rely on large-deviation arguments and need continuity of the map $\mu \mapsto \langle \mu, W \wedge M + (v \circ |\cdot|) \wedge M \rangle$ in the weak topology on the set of probability measures. In the proof of the lower bound, we use an eigenvalue expansion, which needs only local integrability.

(2) The additional assumption that there is a function $\tilde{v}$ with $v \leq \tilde{v}$ and $\int_{B_1(0)} G(0, y)\tilde{v}(|y|)\, dy < \infty$ in the soft-core case is necessary only in our proof of the lower bound. This means that we can handle the case where $v(r) \leq \mathcal{O}(r^{-\varepsilon})$ as $r \downarrow 0$ with $\varepsilon < 2$, but not the case $v(r) \geq Cr^{-\varepsilon}$ with $\varepsilon \in [2, d]$, as the sole requirement $\int_{B_1(0)} \tilde{v}(|y|)\, dy < \infty$ would include.

Let us draw a corollary about the mean of the $N$ Brownian occupation measures,

$$(1.30) \qquad \overline{\mu}_{N,\beta} = \frac{1}{N} \sum_{i=1}^{N} \mu_\beta^{(i)} \in \mathcal{M}_1(\mathbb{R}^d).$$

COROLLARY 1.9.  *Fix* $N \in \mathbb{N}$. *As* $\beta \to \infty$, *the distribution of the mean* $\overline{\mu}_{N,\beta}$ *under* $\widehat{\mathbb{P}}_{N,\beta}^{(\otimes)}$ *satisfies a large deviation principle in the weak topology on* $\mathcal{M}_1(\mathbb{R}^d)$ *with speed* $\beta N$ *and rate function*

$$
\begin{aligned}
(1.31) \qquad I_N^{(\otimes,\mathrm{mean})}(\mu) = \ & \langle W, \mu \rangle + \frac{N}{2} \langle \mu, V\mu \rangle \\
& + \inf_{\mu_1,\dots,\mu_N \in \mathcal{M}_1(\mathbb{R}^d)\,:\,\overline{\mu}=\mu} \Bigg[ \frac{1}{N} \sum_{i=1}^{N} I_1(\mu_i) \\
& \hspace{4em} - \frac{1}{2N} \sum_{i=1}^{N} \langle \mu_i, V\mu_i \rangle \Bigg] - \chi_N^{(\otimes)}.
\end{aligned}
$$

*In particular,* $\overline{\mu}_{N,\beta}$ *converges under* $\widehat{\mathbb{P}}_{N,\beta}^{(\otimes)}$, *as* $\beta \to \infty$, *weakly toward the minimizer(s) of* $I_N^{(\otimes,\mathrm{mean})}$ *on* $\mathcal{M}_1(\mathbb{R}^d)$.

REMARK 1.10.   If one changed the Hartree model by adding all the terms for $i = j$ in the pair-interaction term, then all statements would remain valid after obvious changes in the notation. This model has an additional self-interaction of each path. For this model, stronger statements are possible. Assume that $v$ is such that the map $\mu \mapsto \langle \mu, V\mu \rangle$ is convex on $\mathcal{M}_1(\mathbb{R}^d)$. Then the rate function in Corollary 1.9 (without the term $\frac{1}{2N} \sum_{i=1}^{N} \langle \mu_i, V\mu_i \rangle$) may be identified as $I(\mu) + \langle W, \mu \rangle + \langle \mu, V\mu \rangle$ (minus the normalization) and therefore turns out to be strictly convex.

Let us now comment on the physical relevance of the Hartree model. Recall that the Hartree formula in (1.9) was introduced as an ansatz for studying the ground states of $\mathcal{H}_N$. Theorem 1.7 shows that the Hartree model is the correct positive-temperature model for this ansatz. In this way, its relation to the product ground states of $\mathcal{H}_N$ is analogous to the relation of the canonical ensemble model to its ground state.

The study of the large-$N$ behavior of the canonical ensemble model at *positive* temperature is an important and difficult open problem. We conceive the Hartree model as an ansatz for approaching this task. Indeed, we study this question for the Hartree model in [1]. For this question at *zero* temperature, see Section 1.5.

1.4. *Dirac interaction.* In this section we discuss a further interesting choice of the interaction potential. We restrict now to dimensions $d \in \{2, 3\}$. We do not choose the interaction as a *function*, but as a *measure*, still keeping the singularity at zero. More precisely, we replace $v$ by the Dirac measure at zero, $\delta_0$.

Let us first remark that the canonical ensemble model does not feel this interaction, since $\mathbb{P}(\exists s : B_s^{(i)} = B_s^{(j)}) = 0$ and hence, $\sum_{1 \leq i < j \leq N} \int_0^\beta \delta_0(|B_s^{(i)} - B_s^{(j)}|) \, ds = 0$ almost surely. Hence, this model with $v$ replaced by $\delta_0$ is the same as if $v$ would be replaced by 0.

However, in $d = \{2, 3\}$, the Hartree model is nontrivial and highly interesting. In order to see this, recall the *intersection local time* of the $i$th and the $j$th motion for $i < j$, which is formally defined by

$$(1.32) \qquad \alpha_\beta^{(i,j)}(x) = \frac{1}{\beta^2} \int_0^\beta ds \int_0^\beta dt \, \delta_x(|B_s^{(i)} - B_t^{(j)}|).$$

It is known [17] that there is a continuous stochastic process $(\alpha_\beta^{(i,j)}(x))_{x \in \mathbb{R}^d}$ which justifies this formal definition, that is, for any continuous bounded function $f : \mathbb{R}^d \to \mathbb{R}$, we have the formula

$$\int_{\mathbb{R}^d} f(x) \alpha_\beta^{(i,j)}(x) \, dx = \frac{1}{\beta^2} \int_0^\beta ds \int_0^\beta dt \, f(|B_s^{(i)} - B_t^{(j)}|).$$

Furthermore, $\alpha_\beta^{(i,j)}$ is even continuous in 0, which makes it possible to define the normalized amount of intersection of $B^{(i)}$ and $B^{(j)}$ as the random variable $\alpha_\beta^{(i,j)}(0)$. Hence, we can replace the Hamiltonian $K_{N,\beta}$ in (1.25) by

$$(1.33) \qquad \widetilde{K}_{N,\beta} = \lambda \beta \sum_{1 \leq i < j \leq N} \alpha_\beta^{(i,j)}(0) \qquad \text{for some } \lambda \in (0, \infty).$$

The relevant variational problem should be

$$
\begin{aligned}
\chi_N^{(\delta_0)}(\lambda) = \frac{1}{N} \\
(1.34) \qquad \times \inf_{h_1, \ldots, h_N \in H^1(\mathbb{R}^d) \, : \, \|h_i\|_2 = 1 \, \forall i} \Bigg\{ \sum_{i=1}^N \|\nabla h_i\|_2^2 \\
+ \sum_{i=1}^N \langle W, h_i^2 \rangle + \lambda \sum_{1 \leq i < j \leq N} \langle h_i^2, h_j^2 \rangle \Bigg\}.
\end{aligned}
$$

Using the means of the proof of Lemma 1.3 in Sections 4.1–4.3 below, one can also show the following.

LEMMA 1.11 [Minimizers in (1.34)]. *Fix $\lambda \in (0, \infty)$.*

(i) *The infimum in* (1.34) *is attained.*

(ii) *Any minimizer* $(h_1, \ldots, h_N)$ *of* (1.34) *satisfies the following system of Euler–Lagrange equations:*

$$\Delta h_i = -\lambda_i h_i + W h_i + \lambda h_i \sum_{j \neq i} h_j^2, \qquad i = 1, \ldots, N,$$

*where* $\lambda_i = \|\nabla h_i\|_2^2 + \langle W, h_i^2 \rangle + \lambda \sum_{j \neq i} \langle h_i^2, h_j^2 \rangle$.

(iii) *Let* $(h_1, \ldots, h_N)$ *be a minimizer of* (1.34). *Then every* $h_i$ *is positive almost everywhere. Furthermore, in* $d = 2$, *every* $h_i$ *is continuously differentiable, and the first derivatives are* $\alpha$-*Hölder continuous for any* $\alpha < 1$. *In* $d = 3$, *every* $h_i^p$ *is locally integrable for any* $p \in (0, 3)$.

It is natural to conjecture that, for any $N \in \mathbb{N}$ and any $\lambda \in (0, \infty)$,

$$(1.35) \qquad \lim_{\beta \to \infty} \frac{1}{N\beta} \log \mathbb{E}(e^{-H_{N,\beta} - \widetilde{K}_{N,\beta}}) = -\chi_N^{(\delta_0)}(\lambda),$$

and that the tuple of normalized occupation measures of the motions stands in the analogous relation to the minimizers as in Theorem 1.7. However, we do not know how to prove this conjecture, since the treatment of Brownian intersection local times is technically notoriously difficult. Instead, we offer an analogous result for simple random walks in place of Brownian motions, which we can handle rigorously.

For this purpose, let $(S_t^{(i)})_{t \in [0, \infty)}$ be independent continuous-time simple random walks on $\mathbb{Z}^d$ for $i = 1, \ldots, N$ with generator $-\Delta$ given by $\Delta f(z) = \sum_{y \sim z} (f(y) - f(z))$. For simplicity, we pick some joint initial distribution with compact support.

The normalized intersection local time of the $i$th and the $j$th walk is given as

$$\alpha_\beta^{(i,j)} = \frac{1}{\beta^2} \sum_{z \in \mathbb{Z}^d} \ell_\beta^{(i)}(z) \ell_\beta^{(j)}(z),$$

where $\ell_\beta^{(i)}(z) = \int_0^\beta \mathbb{1}\{S_t^{(i)} = z\} \, dt$ denotes the local times of the $i$th walk up to time $\beta > 0$. In order to have a perfect-scaling property, we restrict to the trap potential $W(x) = |x|^p$ for some $p > d - 2$. We consider the following model:

$$(1.36) \qquad d\widehat{\mathbb{P}}_{N,\beta}^{(\lambda)} = \frac{1}{Z_{N,\beta}^{(\lambda)}} \exp\left(-\frac{1}{\beta} \sum_{i=1}^N \int_0^\beta W(S_t^{(i)}) \, dt \right.$$
$$\left. - \lambda \beta^{(d+p)/(2+p)} \sum_{1 \leq i < j \leq N} \alpha_\beta^{(i,j)} \right) d\mathbb{P}.$$

As usual, $Z_{N,\beta}^{(\lambda)} > 0$ denotes the constant that makes $\widehat{\mathbb{P}}_{N,\beta}^{(\lambda)}$ a probability measure. Note the factor of $1/\beta$ in front of the trap term and the factor of $\beta^{(d+p)/(2+p)}$ in front of the interaction term. The first ensures that the properly rescaled random walks approach some Brownian motions, which makes the model asymptotically equal to the above Brownian model. In order to formulate our result on the large-$\beta$ asymptotic of $\widehat{\mathbb{P}}_{N,\beta}^{(\lambda)}$, we need to introduce the following normalized and rescaled version of $\ell_\beta^{(i)}$:

$$(1.37) \qquad L_\beta^{(i)}(x) = \frac{\xi_\beta^d}{\beta} \ell_\beta^{(i)}(\lfloor x\xi_\beta \rfloor), \qquad x \in \mathbb{R}^d, \text{ where } \xi_\beta = \beta^{1/(2+p)}.$$

Note that $L_\beta^{(i)}$ is a (random) probability density on $\mathbb{R}^d$. Let $\mu_\beta^{(i)}(dx) = L_\beta^{(i)}(x)\,dx$ be the corresponding measure.

**THEOREM 1.12** (Discrete Dirac-interaction model at late times). *Fix $d \geq 2$ and $p > d-2$. Furthermore, let $N \in \mathbb{N}$ and $\lambda \in (0,\infty)$. Then*

(i)

$$(1.38) \qquad \lim_{\beta \to \infty} \beta^{-p/(2+p)} \log Z_{N,\beta}^{(\lambda)} = -N\chi_N^{(\delta_0)}(\lambda).$$

(ii) *As $\beta \to \infty$, the distribution of the tuple $(\mu_\beta^{(1)}, \ldots, \mu_\beta^{(N)})$ of normalized rescaled local times under $\widehat{\mathbb{P}}_{N,\beta}^{(\lambda)}$ on $\mathcal{M}_1(\mathbb{R}^d)^N$ satisfies a large deviation principle with speed $\beta^{p/(2+p)}$ and rate function*

$$
\begin{aligned}
I_N^{(\lambda)}(\mu_1, \ldots, \mu_N) &= \sum_{i=1}^N I_1(\mu_i) + \langle \mathfrak{W}, \mu^\otimes \rangle \\
(1.39) \qquad &\quad + \lambda \sum_{1 \leq i < j \leq N} \left\langle \frac{d\mu_i}{dx}, \frac{d\mu_j}{dx} \right\rangle - N\chi_N^{(\delta_0)}(\lambda),
\end{aligned}
$$

$$\mu_1, \ldots, \mu_N \in \mathcal{M}_1(\mathbb{R}^d),$$

*where $I_1$ is defined in (1.19), $\mu^\otimes = \mu_1 \otimes \cdots \otimes \mu_N$ is the product measure, and we define $I_N^{(\lambda)}(\mu_1, \ldots, \mu_N) = +\infty$ if any of the measures $\mu_1, \ldots, \mu_N$ fails to have a Lebesgue density.*

(iii) *The distribution of $(\mu_\beta^{(1)}, \ldots, \mu_\beta^{(N)})$ under $\widehat{\mathbb{P}}_{N,\beta}^{(\lambda)}$ is attracted by the set of minimizers in (1.34).*

For the proof see Section 3.

Note that Theorem 1.12 is true in any dimension $d \geq 2$, while the Brownian version is well defined only in $d \in \{2, 3\}$.

The choice of the factor $1/\beta$ in front of the trap term in (1.36) is highhanded. Because of this term, one has to assume that $p > d-2$ in order that an appropriate large deviation principle be applicable (see Lemma 1.16). If $1/\beta$ would be replaced by $\beta^{-2/d}$, then the assumption $p > 0$ would suffice.

1.5. *Large-N behavior of the product ground states.* In this section we study our main variational formulas, $\chi_N$ and $\chi_N^{(\otimes)}$, and their minimizers in the limit for diverging number $N$ of particles. In particular, we point out some significant differences between $\chi_N$ and its product state version $\chi_N^{(\otimes)}$ in the soft-core and the hard-core case, respectively.

First we report on recent results by Lieb, Seiringer and Yngvason on the large-$N$ behavior of $\chi_N$. Let the pair functional $v$ be as in (1.2) and assume additionally that $v \geq 0$ and $v(0) > 0$.

We shall replace $v$ by the rescaling $v_N(\cdot) = \xi_N^{-2} v(\cdot \xi_N^{-1})$, for some appropriate $\xi_N$ tending to zero sufficiently fast. Hence, the reach of the repulsion is of order $\xi_N$, and its strength of order $\xi_N^{-2}$. Furthermore, the *scattering length* of $v$, $\alpha(v)$, is rescaled such that $\alpha(v_N) = \alpha(v)\beta_N$ (see Section 1.6.1 below for the definition of the scattering length and some of its properties). If $\beta_N \downarrow 0$ sufficiently fast, this rescaling makes the system dilute, in the sense that $\alpha(v_N) \ll N^{-1/d}$. This means that the interparticle distance is much bigger than the range of the interaction potential strength. More precisely, the decay of $\beta_N$ will be chosen in such a way that the pair-interaction has the same order as the kinetic term. In $d = 3$, this choice of rescaling is motivated by famous experiments for the derivation of Bose–Einstein condensation of a large, but finite, dilute trapped system of $N$ real particles ($^{87}$Rb [4], $^{7}$Li [6], $^{23}$Na [10], but also more recently $^{85}$Rb, $^{41}$K, $^{133}$Cs, hydrogen, metastable triplet $^{4}$He, $^{174}$Yb, $^{85}$Rb$_2$, and $^{6}$Li$_2$). Here the scattering length is of the order $10^{-3}$, whereas $N$ varies from $10^3$ to $10^7$.

The mathematical description of the large-$N$ behavior of $\chi_N$ in this scaling, and hence, the theoretical foundation of the above mentioned physical experiments, has been successfully accomplished in a recent series of papers [22, 24, 25, 26]. It turned out that the well-known *Gross–Pitaevskii formula* adequately describes the limit of the ground states and its energy. This variational formula was first introduced in [19] and [20] and independently in [27] for the study of superfluid Helium. After its importance for the description of Bose–Einstein condensation of dilute gases in magnetic traps was realized, the interest in this formula considerably increased; see [9] for a summary and the monograph [28] for a comprehensive account on Bose–Einstein condensation.

The Gross–Pitaevskii formula has a parameter $\alpha > 0$ and is defined as follows:

$$(1.40) \qquad \chi_\alpha^{(\mathrm{GP})} = \inf_{\phi \in H^1(\mathbb{R}^d):\, \|\phi\|_2 = 1} [\|\nabla \phi\|_2^2 + \langle W, \phi^2 \rangle + 4\pi\alpha \|\phi\|_4^4].$$

It is known [22] that $\chi_\alpha^{(\mathrm{GP})}$ possesses a unique minimizer $\phi_\alpha^{(\mathrm{GP})}$, which is positive and continuously differentiable with Hölder continuous derivatives of order one.

Since $v(0) > 0$, its scattering length $\alpha(v)$ is positive (see Section 1.6.1 below). The condition $\int_{a+1}^{\infty} v(r) r^{d-1} \, dr < \infty$ implies that $\alpha(v) < \infty$. Furthermore, note that the rescaled potential $\xi^{-2} v(\cdot \xi^{-1})$ has scattering length $\xi \alpha(v)$ for any $\xi > 0$.

THEOREM 1.13 (Large-$N$ asymptotic of $\chi_N$ in $d \in \{2,3\}$, [22, 24, 26]). *Assume that $d \in \{2,3\}$, that $v \geq 0$ with $v(0) > 0$, and $\int_{a+1}^{\infty} v(r) r^{d-1} \, dr < \infty$. Replace $v$ by $v_N(\cdot) = \xi_N^{-2} v(\cdot \xi_N^{-1})$ with $\xi_N = 1/N$ in $d = 3$ and $\xi_N^2 = \alpha(v)^{-2} \times e^{-N/\alpha(v)} N \|\phi_{\alpha(v)}^{(\mathrm{GP})}\|_4^{-4}$ in $d = 2$. Let $h_N \in H^1(\mathbb{R}^{dN})$ be the unique minimizer on the right-hand side of (1.6), and define $\phi_N^2 \in H^1(\mathbb{R}^d)$ as the normalized first marginal of $h_N^2$, that is,*

$$\phi_N^2(x) = \int_{\mathbb{R}^{d(N-1)}} h_N^2(x, x_2, \ldots, x_N) \, dx_2 \cdots dx_N, \qquad x \in \mathbb{R}^d.$$

*Then we have*

$$\lim_{N \to \infty} \chi_N = \chi_{\alpha(v)}^{(\mathrm{GP})} \quad and \quad \phi_N^2 \to (\phi_{\alpha(v)}^{(\mathrm{GP})})^2 \qquad in \ weak \ L^1(\mathbb{R}^d)\text{-}sense.$$

In particular, the proofs show that the ground state, $h_N$, approaches the ground state $(\phi_{\alpha(v)}^{(\mathrm{GP})})^{\otimes N}$ if $N$ gets large. In order to obtain the Gross–Pitaevskii formula as the limit of $\chi_N$ also in $d = 2$, the rescaling of $v$ in Theorem 1.13 has to be chosen in such a way that the repulsion strength is the inverse square of the repulsion reach and such that this reach decays exponentially, which is rather unphysical.

In the present paper, we prove the analogue of Theorem 1.13 for the Hartree model in the soft-core case. It turns out that $\chi_N^{(\otimes)}$ in (1.11) also converges toward the Gross–Pitaevskii formula. However, in $d = 2$, it turns out that the potential $v$ has to be rescaled differently. Furthermore, in $d \in \{2,3\}$, the scattering length $\alpha(v)$ is replaced by the number

$$(1.41) \qquad \widetilde{\alpha}(v) := \frac{1}{8\pi} \int_{\mathbb{R}^d} v(|y|) \, dy.$$

THEOREM 1.14 (Large-$N$ asymptotic of $\chi_N^{(\otimes)}$, soft-core case). *Let $d \in \{2,3\}$. Assume that $v$ is a soft-core pair potential with $v \geq 0$ and $v(0) > 0$ and $\widetilde{\alpha}(v) < \infty$. In dimension $d = 3$, additionally assume that (1.13) holds. Replace $v$ by $v_N(\cdot) = N^{d-1} v(N \cdot)$ and let $(h_1^{(N)}, \ldots, h_N^{(N)})$ be any minimizer on the right-hand side of (1.11). Define $\phi_N^2 = \frac{1}{N} \sum_{i=1}^N (h_i^{(N)})^2$. Then we have*

$$\lim_{N \to \infty} \chi_N^{(\otimes)} = \chi_{\widetilde{\alpha}(v)}^{(\mathrm{GP})} \quad and \quad \phi_N^2 \to (\phi_{\widetilde{\alpha}(v)}^{(\mathrm{GP})})^2,$$

*where the convergence of $\phi_N^2$ is in the weak $L^1(\mathbb{R}^d)$-sense and weakly for the probability measures $\phi_n^2(x) \, dx$ toward the measure $(\phi_{\widetilde{\alpha}(v)}^{(\mathrm{GP})})^2(x) \, dx$.*

For the proof see Section 4.5.

Note that, in $d = 3$, the interaction potential is rescaled in the same way in Theorems 1.13 and 1.14. However, the two relevant parameters depend on different properties of the potential (the scattering length, resp. the integral) and have different values, since $\alpha(v) < \widetilde{\alpha}(v)$ (see [26] and Section 1.6.1 below). In particular, for $N$ large enough, the ground state of $\chi_N$ is *not* a product state. This implies the strictness of the inequality in (1.10), for $v$ replaced by $v_N(\cdot) = N^2 v(\cdot N)$. The phenomenon that (unrestricted) ground states are linked with the scattering length has been theoretically predicted for more general $N$-body problems [14], Chapter 14, [29]. Indeed, Landau combined a diagrammatic method (a Born approximation of the scattering length) with Bogoliubov's approximations to almost reconstruct the scattering length from the $L^1$-norm of $v \circ |\cdot|$ in the (nondilute) ground state. However, the relation between the $L^1$-norm and the product ground states was not rigorously known before.

In $d = 2$, a more substantial difference between the large-$N$ behaviors of $\chi_N$ and $\chi_N^{(\otimes)}$ is apparent. Not only the asymptotic relation between the reach and the strength of the repulsion is different, but also the order of this rescaling in dependence on $N$. We can offer no intuitive explanation for this.

Interestingly, in the hard-core case, $\chi_N^{(\otimes)}$ shows a rather different large-$N$ behavior, which we want to roughly indicate in a special case. Assume that $W$ and $v$ are purely hard-core potentials, for definiteness, we take $W = \infty \mathbb{1}_{B_1(0)^c}$ and $v = \infty \mathbb{1}_{[0,a]}$. We replace $v$ by $v_N(\cdot) = v(\cdot/\xi_N)$ for some $\xi_N \downarrow 0$ (a pre-factor plays no role). Then $\chi_N^{(\otimes)}$ is equal to $\frac{1}{N}$ times the minimum over the sum of the principal Dirichlet eigenvalues of $-\Delta$ in $N$ subsets of the unit ball having distance $\geq a\beta_N$ to each other, where the minimum is taken over the $N$ sets. It is clear that the volumes of these $N$ sets should be of order $\frac{1}{N}$, independently of the choice of $\xi_N$. Then their eigenvalues are at least of order $N^{2/d}$. Hence, one arrives at the statement $\liminf_{N\to\infty} N^{-2/d} \chi_N^{(\otimes)} > 0$, that is, $\chi_N^{(\otimes)}$ tends to $\infty$ at least like $N^{2/d}$.

## 1.6. *Preliminaries.*

### 1.6.1. *The scattering length.*

Let us briefly introduce the scattering length of the pair potential, $v$, and its most important properties. For a detailed overview, see [26]. First we turn to $d \geq 3$. Let $u: [0, \infty) \to [0, \infty)$ be a solution of the *scattering equation*,

$$(1.42) \qquad u'' = \tfrac{1}{2} uv \qquad \text{on } (0, \infty), \qquad u(0) = 0.$$

Then the scattering length $\alpha(v) \in [0, \infty]$ of $v$ is defined as

$$(1.43) \qquad \alpha(v) = \lim_{r \to \infty} \left[ r - \frac{u(r)}{u'(r)} \right].$$

If $v(0) > 0$, then $\alpha(v) > 0$, and if $\int_{a+1}^{\infty} v(r) r^{d-1} \, dr < \infty$, then $\alpha(v) < \infty$. In the pure hard-core case, that is, $v = \infty \mathbb{1}_{[0,a)}$, we have $\alpha(v) = a$. It is easily seen from the definition that the scattering length of the rescaled potential $\xi^{-2} v(\cdot \, \xi^{-1})$ is equal to $\xi \alpha(v)$, for any $\xi > 0$.

There is some ambiguity of the choice of $u$ in (1.42); positive multiples of $u$ are also solutions, but the factor drops out in (1.43). We like to normalize $u$ by requiring that $\lim_{R \to \infty} u'(R) = 1$. It is easily seen that (where $\omega_d$ denotes the area of the unit sphere in $\mathbb{R}^d$),

$$
\begin{aligned}
\int_{\mathbb{R}^d} v(|x|) \frac{u(|x|)}{|x|^{d-2}} \, dx &= \omega_d \int_0^{\infty} v(r) u(r) r \, dr \\
&= 2\omega_d \int_0^{\infty} u''(r) r \, dr \\
&= 2\omega_d \lim_{R \to \infty} \left( u'(r) r \big|_0^R - \int_0^R u'(r) \, dr \right) \\
&= 2\omega_d \lim_{R \to \infty} (u'(R) R - u(R)) = 2\omega_d \alpha(v).
\end{aligned}
\tag{1.44}
$$

As a consequence, in dimension $d = 3$, we have $\alpha(v) < \widetilde{\alpha}(v)$. Indeed, $u$ is a nonnegative convex function whose slope is always below one because of $\lim_{R \to \infty} u'(R) = 1$. By $u(0) = 0$, we have that $u(r) < r = r^{d-2}$ for any $r > 0$. With the help of (1.44), we therefore get $8\pi \alpha(v) = 2\omega_d \alpha(v) < \int_{\mathbb{R}^d} v(|x|) \, dx = 8\pi \widetilde{\alpha}(v)$.

In $d = 2$, the definition of the scattering length is slightly different. We treat first the case that $\mathrm{supp}(v) \subset [0, R_*]$ for some $R_* > 0$ and consider, for some $R > R_*$, the solution $u \colon [0, R] \to [0, \infty)$ of the scattering equation

$$
u'' = \tfrac{1}{2} uv \qquad \text{on } [0, R], \qquad u(R) = 1, u(0) = 0.
$$

Then $u(r) = \log \frac{r}{\alpha(v)} / \log \frac{R}{\alpha(v)}$ for $R_* < r < R$ for some $\alpha(v) \geq 0$, which is by definition the scattering length of $v$ in the case that $\mathrm{supp}(v) \subset [0, R_*]$. Note that $\alpha(v)$ does not depend on $R$. Hence,

$$
\log \alpha(v) = \frac{\log r - u(r) \log R}{1 - u(R)}, \qquad R_* < r < R.
$$

For general $v$ (i.e., not necessarily having finite support), $v$ is approximated by compactly supported potentials, and the scattering length of $v$ is put equal to the limit of the scattering lengths of the approximations.

### 1.6.2. Large deviations principles.

For the convenience of our reader, we repeat the notion of a large deviation principle. A family $(X_\beta)_{\beta > 0}$ of random variables $X_\beta$, taking values in a topological vector space $\mathcal{X}$, satisfies the *large deviation upper bound* with *speed* $a_\beta$, where $a_\beta \to \infty$ for $\beta \to \infty$, and rate function $I \colon \mathcal{X} \to [0, \infty]$ if, for any closed subset $F$ of $\mathcal{X}$,

$$
\limsup_{\beta \to \infty} \frac{1}{a_\beta} \log \mathbb{P}(X_\beta \in F) \leq -\inf_{x \in F} I(x),
$$

and it satisfies the *large deviation lower bound* if, for any open subset $G$ of $\mathcal{X}$,

$$\liminf_{\beta \to \infty} \frac{1}{a_\beta} \log \mathbb{P}(X_\beta \in G) \geq -\inf_{x \in G} I(x).$$

If both, upper and lower bound, are satisfied and, in addition, the level sets $\{I \leq c\}$ are compact for any $c \in \mathbb{R}$, then one says that $(X_\beta)_\beta$ satisfies a *large deviation principle*. This notion easily extends to the situation where the distribution of $X_\beta$ is not normalized, but a sub-probability distribution only.

In the proofs of Theorem 1.7 we shall rely on the following principles for the normalized Brownian occupation measures, that is, for certain $\mathcal{M}_1(\mathbb{R}^d)$-valued random variables. For any measurable subset $A$ of $\mathbb{R}^d$, we conceive $\mathcal{M}_1(A)$ as a closed convex subset of the space $\mathcal{M}(A)$ of all finite signed Borel measures on $A$, which is a topological Hausdorff vector space whose topology is induced by the set $\mathcal{C}_b(A)$ of all continuous bounded functions $A \to \mathbb{R}$. Here $\mathcal{C}_b(A)$ is the topological dual of $\mathcal{M}(A)$. The set $\mathcal{M}_1(\mathbb{R}^d)$ inherits this topology from $\mathcal{M}(A)$. When we speak of a large deviation principle of $\mathcal{M}_1(A)$-valued random variables, then we mean a principle on $\mathcal{M}(A)$ with a rate function that is tacitly extended from $\mathcal{M}_1(A)$ to $\mathcal{M}(A)$ with the value $+\infty$.

One of the principles we are going to present is for the Brownian motion in a given bounded set, and the other one for a periodized version of the motion, which we introduce now. Let $R > 0$ and $\Lambda_R = [-R, R]^d \subset \mathbb{R}^d$. For any probability measure $\mu$ on $\mathbb{R}^d$, we denote by $\mu_R \in \mathcal{M}_1(\Lambda_R)$ the periodized version of $\mu$, that is,

$$(1.45) \qquad \mu_R(A) = \mu\left( \bigcup_{k \in \mathbb{Z}^d} (A + 2kR) \right), \qquad A \subset \Lambda_R \text{ measurable.}$$

Note that the shifted cubes $\Lambda_R + 2kR$ with $k \in \mathbb{Z}^d$ are disjoint up to their boundaries, and that their union covers $\mathbb{R}^d$. Recall the Donsker–Varadhan rate function from (1.19). The periodized version of the rate function $I_1$ on $\mathcal{M}_1(\Lambda_R)$ is denoted $I_1^{(R)}$, that is,

$$(1.46) \quad I_1^{(R)}(\mu) = \inf\{I_1(\nu) : \mu = \nu_R \text{ is the periodized version of } \nu\}.$$

LEMMA 1.15 (Large deviations principle for occupation measures [16]). *Fix $d \in \mathbb{N}$. Let $(B_t)_{t \geq 0}$ be the Brownian motion on $\mathbb{R}^d$ with generator $-\Delta$, and let $\mu_\beta(dx) = \frac{1}{\beta} \int_0^\beta \delta_{B_s}(dx)\, ds$ be the normalized occupation measure up to time $\beta > 0$. Fix $R > 0$.*

(i) *The family $(\mu_{\beta,R})_{\beta > 0}$ of $\Lambda_R$-periodizations of $\mu_\beta$ satisfies the large deviations upper bound with speed $\beta$ on $\mathcal{M}_1(\Lambda_R)$ with rate function $I_1^{(R)}$.*

(ii) *For any open bounded set $A \subset \mathbb{R}^d$, the family $(\mu_\beta)_{\beta>0}$ satisfies, under the sub-probability measures $\mathbb{P}(\cdot \cap \{\mathrm{supp}(\mu_\beta) \subset A\})$, the large deviations lower bound with speed $\beta$ on $\mathcal{M}_1(A)$; the corresponding rate function is the restriction of $I$ to the set of probability measures whose support lies in $A$.*

In the proof of Theorem 1.12 in Section 3 we shall rely on related principles for the normalized and rescaled local times $L_\beta^{(i)}$ defined in (1.37), more precisely, for the corresponding measures $\mu_\beta^{(i)}(dx) = L_\beta^{(i)}(x)\,dx$.

LEMMA 1.16 (Large deviations principle for rescaled local times of random walks). *Fix $d \in \mathbb{N}$. Let $(S_t)_{t \in [0,\infty)}$ be a simple random walk on $\mathbb{Z}^d$ with generator $-\Delta$ (the discrete version of the Laplace operator) and local times $\ell_\beta(z) = \int_0^\beta \mathbb{1}\{S_t = z\}\,dt$, and let $1 \ll \xi_\beta \ll \beta^{1/d}$ as $\beta \to \infty$ be some scale function. In $d = 1$ assume that $\xi_\beta \ll \sqrt{\beta}$, in $d = 2$ assume that $\xi_\beta^2 \ll \beta/\log\beta$. Define the normalized and rescaled occupation measure by*

$$\mu_\beta(dx) = \frac{\xi_\beta^d}{\beta} \ell(\lfloor x\xi_\beta \rfloor)\,dx.$$

(i) *The family $(\mu_{\beta,R})_{\beta>0}$ of $\Lambda_R$-periodizations of $\mu_\beta$ satisfies the large deviations upper bound with speed $\beta\xi_\beta^{-2}$ on $\mathcal{M}_1(\Lambda_R)$ with rate function $I_1^{(R)}$.*

(ii) *For any open bounded set $A \subset \mathbb{R}^d$, the family $(\mu_\beta)_{\beta>0}$ satisfies, under the sub-probability measures $\mathbb{P}(\cdot \cap \{\mathrm{supp}(\mu_\beta) \subset A\})$, the large deviations lower bound with speed $\beta\xi_\beta^{-2}$ on $\mathcal{M}_1(A)$; the corresponding rate function is the restriction of $I$ to the set of probability measures whose support lies in $A$.*

PROOF. See [15], Lemma 3.2, for the case of a discrete-time random walk. The proof for the continuous-time setting is very similar and is therefore omitted. □

## 2. Large deviations for the Hartree model: Proof of Theorem 1.7.

In this section we prove Theorem 1.7. We shall proceed according to the well-known Gärtner–Ellis theorem. Therefore, we have to establish the existence of the logarithmic moment generating function of $(\mu_\beta^{(1)}, \ldots, \mu_\beta^{(N)})$, that is, the existence of

$$(2.1) \qquad \Lambda_N^{(\otimes)}(\Phi) = \lim_{\beta \to \infty} \frac{1}{\beta} \log \mathbb{E}_{N,\beta}^{(K)}[\exp(\beta\Phi(\mu_\beta^{(1)}, \ldots, \mu_\beta^{(N)}))],$$

for any element $\Phi$ of the dual of the vector space $\mathcal{M}(\mathbb{R}^d)^N$. Note that every linear continuous functional on $\mathcal{M}_1(\mathbb{R}^d)^N$ is of the form $(\mu_1, \ldots, \mu_N) \mapsto$ $\sum_{i=1}^{N} \langle f_i, \mu_i \rangle$ with $f_1, \ldots, f_N \in \mathcal{C}_b(\mathbb{R}^d)$, the set of bounded continuous functions on $\mathbb{R}^d$.

The main step in the proof of Theorem 1.7 is the following.

PROPOSITION 2.1 (Asymptotic for the cumulant generating function). *For any* $f_1, \ldots, f_N \in \mathcal{C}_b(\mathbb{R}^d)$

$$(2.2) \qquad \lim_{\beta \to \infty} \frac{1}{\beta} \log \mathbb{E}[e^{-H_{N,\beta} - K_{N,\beta}} e^{\beta \langle f, \mu_\beta^\otimes \rangle}] = -N\chi_N^{(\otimes)}(f),$$

*where*

$$(2.3) \qquad \chi_N^{(\otimes)}(f) = \frac{1}{N} \inf_{\mu_1, \ldots, \mu_N \in \mathcal{M}_1(\mathbb{R}^d)} \left[ \sum_{i=1}^{N} I_1(\mu_i) + \langle \mathfrak{W} + \mathfrak{v} - f, \mu^\otimes \rangle \right],$$

*and we wrote* $\mu^\otimes = \mu_1 \otimes \cdots \otimes \mu_N$ *and* $f = f_1 \oplus \cdots \oplus f_N$.

The next three subsections are devoted to the proof of the upper and lower bound in (2.2), respectively. In Section 2.4 we finish the proof of Theorem 1.7.

An outline of our proof is the following. Recall (1.26) and (1.27) to see that

$$(2.4) \qquad \begin{aligned} &\mathbb{E}[e^{-H_{N,\beta} - K_{N,\beta}} e^{\beta \langle f, \mu_\beta^\otimes \rangle}] \\ &= \mathbb{E}\left[ \exp\left\{ -\beta \sum_{i=1}^{N} \langle W, \mu_\beta^{(i)} \rangle - \beta \sum_{1 \le i < j \le N} \langle \mu_\beta^{(i)}, V\mu_\beta^{(j)} \rangle \right. \right. \\ &\qquad\qquad\qquad\qquad\qquad\qquad \left. \left. + \beta \sum_{i=1}^{N} \langle f_i, \mu_\beta^{(i)} \rangle \right\} \right] \\ &= \mathbb{E}[\exp\{ -\beta \langle \mathfrak{W} + \mathfrak{v} - f, \mu_\beta^\otimes \rangle \}]. \end{aligned}$$

We intend to apply the large deviation principle in Lemma 1.15 and Varadhan's lemma, which would immediately yield the result in (2.2). However, there are some technical obstacles to be removed. Because $W$ explodes at infinity and $v$ has a singularity at the origin, both functionals $H_{N,\beta}$ and $K_{N,\beta}$ are not continuous and not bounded in the weak topology. Furthermore, we cannot apply immediately the large deviation principle of Lemma 1.15, because the occupation measures $\mu_\beta^{(i)}$ are not restricted to any bounded domain in $\mathbb{R}^d$. As it concerns the proof of the upper bound, these technical obstacles will be removed in Section 2.1 via a well-known cutting and periodization procedure using the large deviations principle in Lemma 1.15(i). An analogous technique works for the proof of the lower bound in the hard-core case, using Lemma 1.15(ii). However, for proving the lower bound in the soft-core case, we did not succeed in making the principle in Lemma 1.15(ii) applicable. The main reason is the singularity of $v$ at zero, which seems to destroy all necessary semicontinuity properties. Instead, we employ an eigenvalue expansion technique for the $N$ iterated expectations w.r.t. the $N$ motions. Here the additional integrability property of $v$ is necessary.

2.1. *Proof of the upper bound in* (2.2). We consider the large closed box $\Lambda_R = [-R, R]^d$. We divide the probability space into the part on which each motion spends more than $(1-\eta)\beta$ time units in $\Lambda_R$ up to time $\beta$ (the main part) and the remaining part, where this is not satisfied (this part will turn out to be negligible). On the first part, we shall replace each $\mu_\beta^{(i)}$ by its $\Lambda_R$-periodized version $\mu_{\beta,R}^{(i)}$ and control the error. Then $H_{N,\beta}$ turns out to be a continuous and bounded functional of the periodized versions. Also, we replace the functions $W$ and $v$ by their cut-off versions $W_M = W \wedge M$ and $v_M = v \wedge M$, respectively, where $M > 0$ is large. This will enable us to apply Varadhan's lemma. Finally, we let the auxiliary parameters $R$, $\eta$ and $M$ tend to infinity respectively to 0.

We turn to the details. Let auxiliary parameters $\eta > 0$ be small and $R, M > 0$ be large. The expectation on the left-hand side of (2.2) is split into two parts and get the estimation

$$
\begin{aligned}
(2.5) \quad & \mathbb{E}(e^{-H_{N,\beta} - K_{N,\beta} + \beta\langle f, \mu_\beta^\otimes \rangle}) \\
& \leq \mathbb{E}\left(e^{-H_{N,\beta} - K_{N,\beta} + \beta\langle f, \mu_\beta^\otimes \rangle} \prod_{i=1}^N \mathbb{1}\{\mu_\beta^{(i)}(\Lambda_R) > 1 - \eta\}\right) \\
& \quad + N\mathbb{E}(e^{-H_{N,\beta}} \mathbb{1}\{\mu_\beta^{(1)}(\Lambda_R) \leq 1 - \eta\}) e^{-\beta[(N^2/2)\inf v - C_f]},
\end{aligned}
$$

where $C_f = \sum_{i=1}^N \|f_i\|_\infty$. The second term is easily estimated, using that $W \geq 0$. Indeed, we have, on the event $\{\mu_\beta^{(1)}(\Lambda_R) \leq 1 - \eta\}$,

$$
H_{N,\beta} \geq \beta \int_{\Lambda_R^c} W(x_1) \, \mu_\beta^{(1)}(dx_1) \geq \beta\eta \inf_{\Lambda_R^c} W,
$$

and therefore, the second term on the right-hand side of (2.5) is not bigger than $Ne^{-\beta(\eta \inf_{\Lambda_R^c} W + (N^2/2)\inf v - C_f)}$.

In the first term, we first estimate $K_{N,\beta} = \beta\langle \mathfrak{v}, \mu_\beta^\otimes \rangle \geq \beta\langle \mathfrak{v}_M, \mu_\beta^\otimes \rangle$, where $v_M(y) = v(y) \wedge M$ is the cut-off pair potential, and $\mathfrak{v}_M$ is defined as $\mathfrak{v}$ with $v$ replaced by $v_M$. Analogously, we estimate $H_{N,\beta} = \beta\langle \mathfrak{W}, \mu_\beta^\otimes \rangle \geq \beta\langle \mathfrak{W}_M, \mu_\beta^\otimes \rangle$ with analogous notation. This leads to

$$
\begin{aligned}
(2.6) \quad & \mathbb{E}(e^{-H_{N,\beta} - K_{N,\beta} + \beta\langle f, \mu_\beta^\otimes \rangle}) \\
& \leq Ne^{-\beta(\eta \inf_{\Lambda_R^c} W + (N^2/2)\inf v - C_f)} \\
& \quad + \mathbb{E}\left(e^{-\beta\langle \mathfrak{W}_M, \mu_\beta^\otimes \rangle - \beta\langle \mathfrak{v}_M, \mu_\beta^\otimes \rangle + \beta\langle f, \mu_\beta^\otimes \rangle} \prod_{i=1}^N \mathbb{1}\{\mu_\beta^{(i)}(\Lambda_R) > 1 - \eta\}\right).
\end{aligned}
$$

Now we replace $\mu_\beta^\otimes$ by its periodized version $\mu_{\beta,R}^\otimes = \mu_{\beta,R}^{(1)} \otimes \cdots \otimes \mu_{\beta,R}^{(N)}$. In order to estimate the error, we point out that, for any $i \in \{1, \dots, N\}$,

$$
\begin{aligned}
(2.7) \quad & \int_{\mathbb{R}^d} W_M(x_i) \mu_\beta^{(i)}(dx_i) - \int_{\Lambda_R} W_M(x_i) \mu_{\beta,R}^{(i)}(dx_i) \\
& = \sum_{k \in \mathbb{Z}^d} \int_{\Lambda_R} (W_M(x_i + 2Rk) \mu_\beta^{(i)}(d(x_i + 2Rk)) \\
& \qquad\qquad\qquad - W_M(x_i) \mu_\beta^{(i)}(d(x_i + 2Rk))) \\
& = \sum_{k \in \mathbb{Z}^d \setminus \{0\}} \int_{\Lambda_R} \mu_\beta^{(i)}(d(x_i + 2Rk))(W_M(x_i + 2Rk) - W_M(x_i)) \\
& \geq -\eta M.
\end{aligned}
$$

Analogously, we derive the error estimate $\langle f, \mu_\beta^\otimes \rangle \leq \langle f, \mu_{\beta,R}^\otimes \rangle + \eta C_f$. The replacement error for the second term is estimated in a similar way: for any pair $i, j \in \{1, \dots, N\}$ with $i \neq j$, we have

$$
\begin{aligned}
(2.8) \quad & \int_{\mathbb{R}^d} \int_{\mathbb{R}^d} v_M(|x_i - x_j|) \mu_\beta^{(i)}(dx_i) \mu_\beta^{(j)}(dx_j) \\
& - \int_{\Lambda_R} \int_{\Lambda_R} v_M(|x_i - x_j|) \mu_{\beta,R}^{(i)}(dx_i) \mu_{\beta,R}^{(j)}(dx_j) \\
& = \sum_{k,l \in \mathbb{Z}^d \,:\, k \neq l} \int_{\Lambda_R} \int_{\Lambda_R} [v_M(|x_i - x_j + 2R(k-l)|) - v_M(|x_i - x_j|)] \\
& \qquad\qquad\qquad\qquad \times \mu_\beta^{(i)}(d(x_i + 2Rk)) \mu_\beta^{(j)}(d(x_j + 2Rl)) \\
& \geq (\inf v - M) \sum_{k,l \in \mathbb{Z}^d \,:\, k \neq l} \mu_\beta^{(i)}(\Lambda_R + 2Rk) \mu_\beta^{(j)}(\Lambda_R + 2Rl) \\
& \geq -2\eta(M - \inf v),
\end{aligned}
$$

since $-\infty < \inf v \leq v_M \leq M$ and $\mu_\beta^{(i)}(\Lambda_R^c) \leq \eta$. Summarizing, we obtain from (2.6) that

$$
\begin{aligned}
(2.9) \quad & \mathbb{E}(e^{-H_{N,\beta} - K_{N,\beta} + \beta \langle f, \mu_\beta^\otimes \rangle}) \\
& \leq N e^{-\beta(\eta \inf_{\Lambda_R^c} W + N^2/2 \inf v - C_f)} \\
& \qquad + e^{\beta \eta [NM + N^2(M - \inf v) + C_f]} \mathbb{E}(e^{\beta \langle -\mathfrak{W}_M - \mathfrak{v}_M + f, \mu_{\beta,R}^\otimes \rangle}).
\end{aligned}
$$

We now argue that we have

$$
(2.10) \quad \limsup_{\beta \to \infty} \frac{1}{\beta} \log \mathbb{E}(e^{\beta \langle -\mathfrak{W}_M - \mathfrak{v}_M + f, \mu_{\beta,R}^\otimes \rangle}) \leq -N \chi_N^{(\otimes)}(R, M, f),
$$

where

$$
(2.11)
\begin{aligned}
&\chi_N^{(\otimes)}(R, M, f) \\
&= \frac{1}{N} \inf\bigg\{ \sum_{i=1}^{N} \tfrac{1}{2}\|\nabla_R h_i\|_{R,2}^2 + \langle (\mathfrak{W} + \mathfrak{v}_M - f)|_{\Lambda_R^N}, (h^\otimes)^2 \rangle : \\
&\qquad\qquad h_1, \ldots, h_N \in \mathcal{C}^\infty(\Lambda_R), \|h_i\|_{R,2} = 1\, \forall i = 1, \ldots, N \bigg\},
\end{aligned}
$$

where $\|\cdot\|_{R,2}$ is the norm on $L^2(\Lambda_R)$, $\nabla_R$ is the gradient on the torus $\Lambda_R$ (i.e., having periodic boundary condition), and $h^\otimes = h_1 \otimes \cdots \otimes h_N$. To prove (2.10), we apply the upper bound part of Varadhan's integral lemma; see [11], Lemma 4.3.6. Indeed, according to Lemma 1.15(i), every family $(\mu_{\beta,R}^{(i)})_{\beta>0}$ satisfies the large deviations upper bound as $\beta \to \infty$, and the map

$$
(\mu_1, \ldots, \mu_N) \mapsto \langle -\mathfrak{W}_M - \mathfrak{v}_M + f, \mu_1 \otimes \cdots \otimes \mu_N \rangle
$$

is upper semicontinuous (even continuous) on the set of vectors of probability measures on $\Lambda_R$ in the weak topology. Also, using [11], Example 4.2.7, (2.10) is a direct consequence, noting that

$$
(2.12)
\begin{aligned}
&\chi_N^{(\otimes)}(R, M, f) \\
&= \frac{1}{N} \inf\bigg\{ \sum_{i=1}^{N} I_R^{(\mathrm{per})}(\mu_i) + \langle \mathfrak{W}_M + \mathfrak{v}_M - f, \mu_1 \otimes \cdots \otimes \mu_N \rangle : \\
&\qquad\qquad \mu_1, \ldots, \mu_N \in \mathcal{M}_1(\Lambda_R) \bigg\}.
\end{aligned}
$$

From (2.9) and (2.10) we obtain that, for any $\eta, R, M > 0$,

the l.h.s. of (2.2)

$$
(2.13)
\begin{aligned}
\leq -\min\bigg\{ &\eta \inf_{\Lambda_R^c} W + \frac{N^2}{2} \inf v - C_f, \\
&-\eta[NM + N^2(M - \inf v) + C_f] + N\chi_N^{(\otimes)}(R, M, f) \bigg\}.
\end{aligned}
$$

On the right-hand side of (2.13), we let first $R \to \infty$, then $\eta \downarrow 0$ and finally $M \to \infty$. Recall that $\lim_{R\to\infty} \inf_{\Lambda_R^c} W = \infty$, according to our assumption in (1.1). It is easily seen that the proof of the upper bound in (2.2) is finished as soon as we have shown that, for some $C > 0$ which does not depend on $M$,

$$
(2.14)
\begin{aligned}
&\liminf_{R\to\infty} \chi_N^{(\otimes)}(R, M, f) \geq -\frac{C}{M} + \chi_N^{(\otimes)}(M, f) \quad \text{and} \\
&\liminf_{M\to\infty} \chi_N^{(\otimes)}(M, f) \geq \chi_N^{(\otimes)}(f),
\end{aligned}
$$

where $\chi_N^{(\otimes)}(M, f)$ is defined as $\chi_N^{(\otimes)}(f)$ in (2.3) with $W$ and $v$ replaced by $W_M$ and $v_M$, respectively.

Certainly, we may assume that $\limsup_{M \to \infty} \limsup_{R \to \infty} \chi_N^{(\otimes)}(R, M, f) < \infty$. Introduce

$$(2.15) \qquad \varepsilon_{R,M} = N \frac{\chi_N^{(\otimes)}(R, M, f) + C_f}{\inf_{\Lambda_{R-1}^c} W_M},$$

and note that $\limsup_{R \to \infty} \varepsilon_{R,M} \leq \frac{C}{M}$ for some $C > 0$, not depending on $M$. To prove the first assertion in (2.14), fix $R > 0$ and let $h_1, \ldots, h_N \in \mathcal{C}^2(\Lambda_R)$ be a vector of approximately minimizing functions for the right-hand side of (2.11). We shall construct $L^2$-normalized functions $\widetilde{h}_1, \ldots, \widetilde{h}_N \in H^1(\mathbb{R}^d)$ such that, for some $C > 0$, not depending on $M$ nor on $R$,

$$(2.16) \qquad \begin{aligned} &\sum_{i=1}^{N} \tfrac{1}{2} \|\nabla h_i\|_{R,2}^2 + \langle (\mathfrak{W}_M + \mathfrak{v}_M - f)|_{\Lambda_R^N}, (h^{\otimes})^2 \rangle \\ &\geq -C \varepsilon_{R,M} + \sum_{i=1}^{N} \tfrac{1}{2} \|\nabla \widetilde{h}_i\|_2^2 + \langle \mathfrak{W} + \mathfrak{v}_M - f, (\widetilde{h}^{\otimes})^2 \rangle. \end{aligned}$$

Passing to the infimum over all vectors of $L^2$-normalized functions $\widetilde{h}_1, \ldots, \widetilde{h}_N \in H^1(\mathbb{R}^d)$, and letting $R \to \infty$, we then arrive at the first assertion in (2.14).

In order to show (2.16), pick some $\phi \in \mathcal{C}^\infty(\mathbb{R}^d, [0, 1])$ with $\mathrm{supp}(\phi) \subset \Lambda_R$ and $\phi|_{\Lambda_{R-1}} = \mathbb{1}_{\Lambda_{R-1}}$. Consider the functions $\widetilde{h}_i = h_i \phi / \|h_i \phi\|_{R,2}$, trivially extended to $\mathbb{R}^d$. Hence, $\widetilde{h}_i \in H^1(\mathbb{R}^d)$ for any $i = 1, \ldots, N$. Note that

$$(2.17) \qquad \begin{aligned} \int_{\Lambda_R \setminus \Lambda_{R-1}} h_i(x)^2 \, dx &\leq \frac{\int_{\mathbb{R}^d} W_M(x) h^2(x) \, dx}{\inf_{\Lambda_{R-1}^c} W_M} \\ &\leq N \frac{\chi_N^{(\otimes)}(R, M, f) + C_f}{\inf_{\Lambda_{R-1}^c} W_M} = \varepsilon_{R,M}. \end{aligned}$$

In particular, $\|h_i \phi\|_2^2 \geq 1 - \varepsilon_{R,M}$.

We now show (2.16). Certainly, we may assume that there is a $C > 0$, not depending on $R$ nor on $M$, such that $|\nabla \phi|^2 \leq C$ and $\|\nabla h_i\|_{R,2}^2 < C$ for any $i \in \{1, \ldots, N\}$ and any large $R > 0$ and any sufficiently large $M > 0$. Then we estimate

$$\|\nabla h_i\|_{R,2}^2 - \|\nabla \widetilde{h}_i\|_2^2$$
$$= -\frac{(1 - \|h_i\phi\|_2^2)}{\|h_i\phi\|_2^2}\|\nabla h_i\|_{R,2}^2$$
$$+ \frac{1}{\|h_i\phi\|_2^2}\int_{\Lambda_R\setminus\Lambda_{R-1}}[(1-\phi^2)|\nabla h_i|^2 - h_i^2|\nabla\phi|^2 - 2\phi h_i\nabla h_i\cdot\nabla\phi]$$

$$(2.18) \quad \geq \frac{-\varepsilon_{R,M}}{1-\varepsilon_{R,M}}C$$
$$- \frac{C}{1-\varepsilon_{R,M}}\Big(\int_{\Lambda_R\setminus\Lambda_{R-1}} h_i^2(x)\,dx + \int_{\Lambda_R\setminus\Lambda_{R-1}} |h_i(x)||\nabla h_i(x)|\,dx\Big)$$
$$\geq \frac{-\varepsilon_{R,M}}{1-\varepsilon_{R,M}}C$$
$$- \frac{C}{1-\varepsilon_{R,M}}\Big(\int_{\Lambda_R\setminus\Lambda_{R-1}} h_i^2(x)\,dx + \Big(\int_{\Lambda_R\setminus\Lambda_{R-1}} h_i^2(x)\,dx\Big)^{1/2}\|\nabla h_i\|_2\Big).$$

Now use (2.17) to estimate the right-hand side from below against $-\varepsilon_{R,M}\widetilde{C}$ for some $\widetilde{C}$, not depending on $M$ nor on $R$.

In order to derive a suitable estimate for the other parts of the functionals, note that

$$(2.19) \quad (h^{\otimes})^2 - (\widetilde{h}^{\otimes})^2 = \sum_{l=1}^{N}\Big(\prod_{i=1}^{l-1}\widetilde{h}_i^2\Big)\Big(\prod_{i=l}^{N} h_i^2\Big)\Big[\frac{(1-\phi^2)}{\|h_l\phi\|_2^2} - \frac{(1-\|h_l\phi\|_2^2)}{\|h_l\phi\|_2^2}\Big].$$

For the two terms between square brackets in (2.19), use the bounds

$$0 \leq \frac{1-\phi^2}{\|h_l\phi\|_2^2} \leq \mathbb{1}_{\Lambda_R\setminus\Lambda_{R-1}} \quad \text{and} \quad 0 \leq \frac{1-\|h_l\phi\|^2}{\|h_l\phi\|_2^2} \leq \frac{\varepsilon_{R,M}}{1-\varepsilon_{R,M}},$$

and recall that $\limsup_{M\to\infty}\limsup_{R\to\infty}\chi_N^{(\otimes)}(R,M,f) < \infty$ in order to easily arrive at the estimate in (2.16).

We now prove the second assertion in (2.14). Let $(h_{1,M},\ldots,h_{N,M})$ be minimizers in the definition of $\chi_N^{(\otimes)}(M,f)$. Along suitable subsequences we have, for any $i \in \{1,\ldots,N\}$,

$$h_{i,M} \longrightarrow h_i \qquad \text{as } M\to\infty,$$

for some $h_i \in H^1(\mathbb{R}^d)$ satisfying $\|h_i\|_2 = 1$. The convergence is weak in $L^2(\mathbb{R}^d)$, weak for the gradients in $L^2(\mathbb{R}^d)$, weak in the sense of probability measures and strongly in $L^2(\mathbb{R}^d)$ on compacts. In particular,

$$\liminf_{M\to\infty}(\tfrac{1}{2}\|\nabla h_{i,M}\|_2^2 - \langle f, h_{i,M}^2\rangle) \geq \tfrac{1}{2}\|\nabla h_i\|_2^2 - \langle f, h_i^2\rangle.$$

For any fixed $M' > 0$, we have

$$\liminf_{M\to\infty}\langle h_{i,M}^2, v_M h_{j,M}^2\rangle \geq \liminf_{M\to\infty}\langle h_{i,M}^2, v_{M'} h_{j,M}^2\rangle = \langle h_i^2, v_{M'} h_j^2\rangle.$$

Letting $M' \uparrow \infty$, the monotone convergence theorem implies that

$$\liminf_{M \to \infty} \langle h_{i,M}^2, v_M h_{j,M}^2 \rangle \geq \liminf_{M' \uparrow \infty} \langle h_i^2, v_{M'} h_j^2 \rangle = \langle h_i^2, v h_j^2 \rangle.$$

In the same way, we see that $\liminf_{M \to \infty} \langle W_M, h_{i,M}^2 \rangle \geq \langle W, h_i^2 \rangle$ for any $i$. If we now pass to the infimum over all $(h_1, \ldots, h_N)$, we arrive at the second assertion of (2.14). This finishes the proof of the upper bound in (2.2).

2.2. *Proof of the lower bound in* (2.2), *hard-core case.* We handle first the case when $v$ is a hard-core interaction potential. Recall the parameter $a = \inf\{r > 0 : v(r) < \infty\} \in [0, \infty)$ from (1.2). Fix a family of open bounded sets $A_1, \ldots, A_N \subset \mathbb{R}^d$ whose mutual pairwise distance is bigger than $\eta$, where $\eta > a$ is close to $a$. We shall use the lower bound

$$(2.20) \quad \begin{aligned} &\mathbb{E}(e^{-H_{N,\beta} - K_{N,\beta}} e^{\beta \langle f, \mu_\beta^\otimes \rangle}) \\ &\geq \mathbb{E}\left( e^{-\beta \langle \mathfrak{W}, \mu_\beta^\otimes \rangle - \beta \langle \mathfrak{v}, \mu_\beta^\otimes \rangle} e^{\beta \langle f, \mu_\beta^\otimes \rangle} \prod_{i=1}^N \mathbb{1}\{\operatorname{supp}(\mu_\beta^{(i)}) \subset A_i\} \right). \end{aligned}$$

Note that, on the event $\bigcap_{i=1}^N \{\operatorname{supp}(\mu_\beta^{(i)}) \subset A_i\}$, the map

$$(\mu_1, \ldots, \mu_N) \mapsto \langle \mathfrak{W} + \mathfrak{v} - f, \mu^\otimes \rangle$$

is bounded since the sets $A_1, \ldots, A_N$ are bounded, and since the supports of the measures $\mu_\beta^{(1)}, \ldots, \mu_\beta^{(N)}$ are bounded away from each other by at least $\eta$ [recall from (1.2) that $v$ is bounded on $(\eta, \infty)$]. This map is also bounded and continuous on the set of $(\mu_1, \ldots, \mu_N) \in \mathcal{M}_1(\mathbb{R}^d)^N$ such that $\operatorname{supp}(\mu_i) \subset A_i$ for all $i$. Hence, Varadhan's lemma may be applied and yields, also using [11], Example 4.2.7,

$$(2.21) \quad \begin{aligned} &\liminf_{\beta \to \infty} \frac{1}{\beta} \log[\text{l.h.s. of } (2.20)] \\ &\geq -\inf\left\{ \sum_{i=1}^N I(\mu_i) + \langle \mathfrak{W} + \mathfrak{v} - f, \mu^\otimes \rangle : \mu_1, \ldots, \mu_N \in \mathcal{M}_1(\mathbb{R}^d), \right. \\ &\hspace{6cm} \left. \forall i : \operatorname{supp}(\mu_i) \subset A_i \right\}. \end{aligned}$$

Now we substitute $h_i^2(x)\, dx = \mu_i(dx)$. According to [8], Proposition 3.29, we can restrict to the infimum on $h_i$ in $\mathcal{C}_c^\infty(A_i)$ the set of infinitely often differentiable functions $\mathbb{R}^d \to \mathbb{R}$ having compact support in $A_i$. Also, passing to the infimum over all admissible sets $A_1, \ldots, A_N$, we see that the left-hand side of (2.21) is not smaller than

$$
\begin{aligned}
-\chi_N^{(\otimes)}(f,\eta) \equiv -\frac{1}{N} \inf\Bigg\{ &\sum_{i=1}^N \|\nabla h_i\|_2^2 + \langle \mathfrak{W} + \mathfrak{v} - f, (h^2)^\otimes \rangle: \\
&h_1, \ldots, h_N \in \mathcal{C}_c^\infty(\mathbb{R}^d), \\
&\|h_1\|_2 = \cdots = \|h_N\|_2 = 1, \forall\, i \neq j: \\
&\qquad \mathrm{dist}(\mathrm{supp}(h_i), \mathrm{supp}(h_j)) \geq \eta \Bigg\},
\end{aligned}
\tag{2.22}
$$

where we adapted the notation $(h^2)^\otimes = (h_1)^2 \otimes \cdots \otimes (h_N)^2$. Now we consider the limit as $\eta \downarrow a$. It is clear that

$$
\begin{aligned}
\lim_{\eta \downarrow a} \chi_N^{(\otimes)}(f,\eta) = \frac{1}{N} \inf\Bigg\{ &\sum_{i=1}^N \|\nabla h_i\|_2^2 + \sum_{i=1}^N \Big\langle W + \sum_{j \neq i} V h_j^2 - f_i, h_i^2 \Big\rangle: \\
&h_1, \ldots, h_N \in \mathcal{C}_c^\infty(\mathbb{R}^d), \\
&\|h_1\|_2 = \cdots = \|h_N\|_2 = 1, \forall\, i \neq j: \\
&\qquad \mathrm{dist}(\mathrm{supp}(h_i), \mathrm{supp}(h_j)) > a \Bigg\}.
\end{aligned}
\tag{2.23}
$$

Now it is easy to see that the right-hand side of (2.23) is equal to $\chi_N^{(\otimes)}(f)$. Indeed, let $(\widetilde{h}_1, \ldots, \widetilde{h}_N) \in H^1(\mathbb{R}^d)$ be a minimizer of the variational formula in (2.3). According to Remark 1.4(ii) (which is only for $f_i \equiv 0$, but nevertheless applies also here), the supports of $\widetilde{h}_1, \ldots, \widetilde{h}_N$, which we denote $\Omega_1, \ldots, \Omega_N$, have distance $\geq a$ to each other. Furthermore, in the case $a = 0$, the interiors of $\Omega_1, \ldots, \Omega_N$ are disjoint. Approaching $\widetilde{h}_1, \ldots, \widetilde{h}_N$ with functions whose supports have distances $> a$ to each other, we obtain

$$
\begin{aligned}
\chi_N^{(\otimes)}(f) = \frac{1}{N} \inf\Bigg\{ &\sum_{i=1}^N \|\nabla h_i\|_2^2 + \sum_{i=1}^N \Big\langle W + \sum_{j \neq i} V h_j^2 - f_i, h_i^2 \Big\rangle: \\
&h_i \in \mathcal{C}_c^\infty(\Omega_i) \text{ with } \|h_i\|_2 = 1\, \forall\, i \Bigg\}.
\end{aligned}
$$

Now one sees that this is equal to the right-hand side of (2.23). This ends the proof of the lower bound in (2.2) in the case of a hard-core potential $v$.

2.3. *Proof of the lower bound in* (2.2), *soft-core case.* Now we turn to the proof of the lower bound in (2.2) in the case of a soft-core potential $v$. The proof goes via an iteration of $N$ eigenvalue expansions for the $N$ expectations with respect to the $N$ Brownian motions. Let us first handle the expectation with respect to the $N$th motion, which we denote by $\mathbb{E}^{(N)}$. We only consider those terms that depend on the $N$th motion and work almost surely with respect to the first $N - 1$ motions.

First we estimate the expectation under interest in terms of the same expectation with respect to some Brownian bridge in the time interval $[0, \beta + 1]$ instead of the Brownian motion on $[0, \beta]$. Let $p_1(x, y)$ denote the standard transition density of Brownian motion at time 1 from $x \in \mathbb{R}^d$ to $y \in \mathbb{R}^d$. Note that $1 \geq p_1(x, y)$ for any $x, y \in \mathbb{R}^d$. Denote by $\nu_N$ the initial distribution of the $N$th motion (i.e., the distribution of $B_0^{(N)}$) and by $\mathbb{E}_x^{(N)}$ the expectation with respect to the motion started at $x \in \mathbb{R}^d$. We abbreviate $q_N = W - f_N + \sum_{j<N} V \mu_\beta^{(j)}$. Note that $q_N$ is a random potential which is locally integrable in $\mathbb{R}^d$, almost surely. Fix some $R > 0$. Then we can estimate

$$
\begin{aligned}
&\mathbb{E}^{(N)}[e^{-\beta\langle q_N, \mu_\beta^{(N)} \rangle}] \\
&\quad \geq \int \nu_N(dx) \mathbb{E}_x^{(N)}[e^{-\beta\langle q_N, \mu_\beta^{(N)} \rangle} p_1(B_\beta^{(N)}, x)] \\
(2.24) &\quad \geq e^{-\|f_N\|_\infty + (N-1)\inf v} \int \nu_N(dx) \mathbb{E}_x^{(N)}[e^{-(\beta+1)\langle q_N, \mu_{\beta+1}^{(N)} \rangle} \delta_x(B_{\beta+1}^{(N)})] \\
&\quad \geq e^{-\|f_N\|_\infty + (N-1)\inf v} \\
&\qquad \times \int \nu_N(dx) \mathbb{E}_x^{(N)}[e^{-(\beta+1)\langle q_N, \mu_{\beta+1}^{(N)} \rangle} \mathbb{1}\{\mathrm{supp}(\mu_{\beta+1}^{(N)}) \subset \Lambda_R\} \delta_x(B_{\beta+1}^{(N)})].
\end{aligned}
$$

The right-hand side of (2.24) can be represented in terms of an eigenvalue expansion. Recall that the potential $q_N$ is integrable on $\Lambda_R$. Let $(\lambda_k)_{k \in \mathbb{N}}$ be the sequence of eigenvalues of the operator $\Delta - q_N$ in $\Lambda_R$ with zero boundary condition, and let $(e_k)_{k \in \mathbb{N}}$ be an orthonormal sequence of corresponding eigenfunctions. We assume that $\lambda_1$ is the principal eigenvalue and that $e_1$ is positive in $\Lambda_R$. Then we have from (2.24) that

$$
\begin{aligned}
(2.25) \quad \mathbb{E}^{(N)}[e^{-\beta\langle q_N, \mu_\beta^{(N)} \rangle}] &\geq e^{-\|f_N\|_\infty + (N-1)\inf v} \sum_{k \in \mathbb{N}} e^{(\beta+1)\lambda_k} \langle \nu_N, e_k^2 \rangle \\
&\geq e^{-\|f_N\|_\infty + \inf v} e^{(\beta+1)\lambda_1} \langle \nu_N, e_1^2 \rangle.
\end{aligned}
$$

Fix some bounded $L^2$-normalized function $h_N \in H^1(\mathbb{R}^d)$ satisfying $\mathrm{supp}(h_N) \subset \Lambda_R$, then we may estimate

$$
\begin{aligned}
(2.26) \quad \lambda_1 &\geq -\|\nabla h_N\|_2^2 - \langle q_N, h_N^2 \rangle \\
&\geq -\|\nabla h_N\|_2^2 - \langle W - f_N, h_N^2 \rangle - \|h_N\|_\infty^2 \sum_{j<N} \int_{\Lambda_R} V \mu_\beta^{(j)}(dx) \\
&\geq -C_R,
\end{aligned}
$$

where $C_R$ is nonrandom and depends only on $\sup_{\Lambda_R} W$, $\|f_N\|_\infty$, $h_N$ and $\int_{\Lambda_R} v(|x|)\, dx$. Hence, we obtain from (2.25) that

$$
\begin{aligned}
(2.27) \quad &\mathbb{E}^{(N)}[e^{-\beta\langle q_N, \mu_\beta^{(N)} \rangle}] \\
&\quad \geq e^{-\|f_N\|_\infty + (N-1)\inf v} e^{-C_R} e^{-\beta[\|\nabla h_N\|_2^2 + \langle q_N, h_N^2 \rangle]} \langle \nu_N, e_1^2 \rangle.
\end{aligned}
$$

The technical difficulty is now to find a positive lower bound for $\langle \nu_N, e_1^2 \rangle$ that does not depend on the first $N - 1$ Brownian motions. Assume that $R > 0$ is so large that $\nu_N(\Lambda_{R-1}) > 0$. We now use Harnack's inequality to obtain a pointwise lower bound for $e_1$ in $\Lambda_{R-1}$. We first check that the random potential $q_N$ lies in the Kato class (see [8], Chapter 3, for more on the Kato class). It is clear that we only have to check this for the random potential $V\mu_\beta^{(j)}$ for $j < N$. We now verify that $\lim_{r\downarrow 0} \psi(r) = 0$, where

$$(2.28) \qquad \psi(r) = \sup_{\mu \in \mathcal{M}_1(\mathbb{R}^d)} \sup_{x \in \Lambda_R} \int_{|y-x| < r} G(x,y) V\mu(y)\, dy.$$

This is seen as follows. First note that it suffices to show that

$$\lim_{r\downarrow 0} \sup_{w \in \mathbb{R}^d} \int_{|y| < r} G(0,y) v(|y-w|)\, dy = 0.$$

Recall that we are under the assumption $\int_{B_\varepsilon(0)} G(0,y)\widetilde{v}(|y|)\, dy < \infty$, where $G$ is the Green's function of the free Brownian motion, and $v \le \widetilde{v}$ on $(0,\varepsilon)$, and $\widetilde{v}$ is decreasing. After possible alteration of $\varepsilon$, we can extend $\widetilde{v}$ to a decreasing function $(0,\infty) \to (0,\infty)$ such that $v \le \widetilde{v}$ on $(0,\infty)$. We choose $r \in (0,\varepsilon/3)$. We distinguish the cases $|w| > 2r$ and $|w| \le 2r$. In the first case, we have, for $|y| < r$, that $|y-w| > r > |y|$ and, therefore, $v(|y-w|) \le \widetilde{v}(|y-w|) \le \widetilde{v}(|y|)$. Hence,

$$\sup_{|w|>2r} \int_{|y|<r} G(0,y) v(|y-w|)\, dy \le \int_{|y|<r} G(0,y)\widetilde{v}(|y|)\, dy,$$

and this vanishes as $r \downarrow 0$.

In the case $|w| \le 2r$, we split the integration over $|y| < r$ into the part where $|y| < |w|/2$ and the part where $|y| \ge |w|/2$. On the first part, we have $|y-w| \ge |w| - |y| \ge |y|$ and, therefore, $v(|y-w|) \le \widetilde{v}(|y-w|) \le \widetilde{v}(|y|)$. Therefore, the first part can be estimated by

$$\sup_{|w|\le 2r} \int_{|y|<|w|/2} G(0,y) v(|y-w|)\, dy \le \int_{|y|<r} G(0,y)\widetilde{v}(|y|)\, dy,$$

which vanishes as $r \downarrow 0$. On the area where $|y| < r$ and $|y| \ge |w|/2$, we can estimate $|y-w| \le |y| + |w| \le 3|y|$ and therefore find a constant $C > 0$ such that $G(0,y) \le CG(0,y-w)$, and this means that we can estimate

$$\sup_{|w|\le 2r} \int_{|w|/2<|y|<r} G(0,y) v(|y-w|)\, dy$$

$$\le C \sup_{|w|\le 2r} \int_{|w|/2<|y|<r} G(0,y-w)\widetilde{v}(|y-w|)\, dy$$

$$\le C \int_{|y|<3r} G(0,y)\widetilde{v}(|y|)\, dy,$$

and the proof of $\lim_{r\downarrow 0} \psi(r) = 0$ is finished. In particular, this means that $q_N$ is in the Kato class.

According to [8], Theorem 5.18, there is a constant $A > 0$, only depending on $R$ and $q_N$, such that

$$(2.29) \qquad \inf_{\Lambda_{R-1}} e_1 \geq A \sup_{\Lambda_{R-1}} e_1.$$

A closer inspection of the proofs of [8], Proposition 5.16, Theorems 5.17 and 5.18, shows that the constant $A$ does not really depend on the function $q_N$, but only on the numbers $\sup_{\Lambda_R} W$, $\|f_N\|_\infty$ and on the small-$r$ behavior of the function $\psi$ defined in (2.28). Hence, $A$ is nonrandom and does not depend on the first $N-1$ motions.

So we have $\langle \nu_N, e_1^2 \rangle \geq A^2 \nu_N(\Lambda_{R-1}) \sup_{\Lambda_{R-1}} e_1^2$. Now we get that $\sup_{\Lambda_{R-1}} e_1^2 \geq (2|\Lambda_{R-1}|)^{-1}$. To see this, note that otherwise we would have $\int_{\Lambda_{R-1}^c} e_1^2 \geq \frac{1}{2}$ and, therefore,

$$(2.30) \quad \begin{aligned} -\lambda_1 &\geq -\|f_N\|_\infty + (N-1)\inf v + \langle W, e_1^2 \rangle \\ &\geq -\|f_N\|_\infty + (N-1)\inf v + \tfrac{1}{2}\inf_{\Lambda_{R-1}^c} W \to \infty \qquad \text{as } R \to \infty, \end{aligned}$$

according to the assumption on $W$ in (1.1). But, as we shall see at the end of the proof, this is impossible if $h_N$ is chosen appropriately. Summarizing, there is a constant $C > 0$, only depending on $N$, $W$, $v$, $f_N$, $\nu_N$ and $R$, such that, almost surely with respect to the first $N-1$ motions,

$$(2.31) \quad \begin{aligned} &\mathbb{E}^{(N)}[e^{-\beta\langle q_N, \mu_\beta^{(N)}\rangle}] \\ &\geq C \exp\left\{-\beta\left[\|\nabla h_N\|_2^2 + \left\langle W - f_N + \sum_{j<N} V\mu_\beta^{(j)}, h_N^2 \right\rangle\right]\right\}. \end{aligned}$$

We iterate now this argument for the $i$th motion for $i = N-1, N-2, \ldots, 1$. For doing this, we have to replace the random potential $q_N$ by $q_i = W - f_i + \sum_{j<i} V\mu_\beta^{(j)} + \sum_{j>i} Vh_j^2$ and obtain, for any bounded $L^2$-normalized function $h_i \in H^1(\mathbb{R}^d)$ satisfying $\text{supp}(h_i) \subset \Lambda_R$, almost surely with respect to the first $i-1$ motions,

$$(2.32) \quad \begin{aligned} &\mathbb{E}^{(i)}[e^{-\beta\langle q_i, \mu_\beta^{(i)}\rangle}] \\ &\geq C \exp\left\{-\beta\left[\|\nabla h_i\|_2^2 + \left\langle W - f_i + \sum_{j<i} V\mu_\beta^{(j)} + \sum_{j>i} Vh_j^2, h_i^2 \right\rangle\right]\right\}, \end{aligned}$$

where $C$ does not depend on $\beta$ nor on the motions. This gives, recalling (1.27),

$$\begin{aligned} &\liminf_{\beta\to\infty} \frac{1}{\beta} \log \mathbb{E}[e^{-H_{N,\beta} - K_{N,\beta} + \langle f, \mu_\beta\rangle}] \\ &\geq -\left[\sum_{i=1}^N \|\nabla h_i\|_2^2 + \sum_{i=1}^N \langle W - f_i, h_i^2 \rangle + \sum_{i<j} \langle h_i^2, Vh_j^2 \rangle\right]. \end{aligned}$$

Now maximize the right-hand side over all choices of bounded $L^2$-normalized functions $h_1, \ldots, h_N \in H^1(\mathbb{R}^d)$ satisfying $\mathrm{supp}(h_i) \subset \Lambda_R$. It is easy to see that, in the limit $R \to \infty$, we obtain that the maximum of the right-hand side tends to $-N\chi_N^{(\otimes)}(f)$. Furthermore, one can see that it is possible to choose approximate maximizers $h_1, \ldots, h_N$ (depending on $R$) such that $\limsup_{R \to \infty} [\|\nabla h_i\|_2^2 + \langle W, h_N^2 \rangle] < \infty$ and $\limsup_{R \to \infty} \|h_i\|_\infty < \infty$. In particular, the eigenvalue $\lambda_1$ introduced below (2.24) satisfies $\limsup_{R \to \infty} (-\lambda_1) < \infty$, and this explains why (2.30) is not possible [see (2.26)].

This completes the proof of the lower bound in (2.2) in the case of a soft-core potential $v$.

### 2.4. *Finish of the proof of Theorem* 1.7.

Proposition 2.1 implies the existence of the logarithmic moment generating function $\Lambda_N^{(\otimes)}$ in (2.1) and identifies it with the Legendre transform of the function $I_N^{(\otimes)}$ defined in (1.29). In particular, Theorem 1.7(i) is implied.

Now we prove the large deviations principle in Theorem 1.7(ii). We use the Gärtner–Ellis theorem (see [11], Corollary 4.5.27). For doing this, it suffices to show that the family of tuples $(\mu_\beta^{(1)}, \ldots, \mu_\beta^{(N)})$ is exponentially tight under $\mathbb{P}_{N,\beta}^{(\otimes)}$ as $\beta \to \infty$, and that $\Lambda_N^{(\otimes)}$ is Gâteau-differentiable. The first condition is verified as follows. We follow the technique of the proof of [11], Lemma 6.2.6. Pick sequences $R_k \to \infty$ and $\varepsilon_k \downarrow 0$ and put $K = \bigcap_{k \in \mathbb{N}} \{\mu \in \mathcal{M}_1(\mathbb{R}^d) : \mu(\Lambda_{R_k}^c) \le \varepsilon_k\}$. The Portmanteau theorem implies that $K$ is closed, and Prohorov's theorem implies that $K$ is relatively compact, hence, $K$ is compact. Now fix $i \in \{1, \ldots, N\}$ and note that

$$\mathbb{P}_{N,\beta}^{(\otimes)}(\mu_\beta^{(i)} \notin K) \le \sum_{k \in \mathbb{N}} \mathbb{P}_{N,\beta}^{(\otimes)}(\mu_\beta^{(i)}(\Lambda_{R_k}^c) > \varepsilon_k) \le \sum_{k \in \mathbb{N}} \mathbb{P}_{N,\beta}^{(\otimes)}\left(\langle W, \mu_\beta^{(i)} \rangle \ge \varepsilon_k \inf_{\Lambda_{R_k}} W\right).$$

Now let a large $R > 0$ be given. We additionally require that $\varepsilon_k \inf_{\Lambda_{R_k}} W \ge Rk$ for all $k \in \mathbb{N}$ [here we use our assumption in (1.1)]. We denote $K$ now by $K_R$. Using the fact that $W \ge 0$ and that $v$ is bounded from below, it is easy to derive the existence of some constant $C \in \mathbb{R}$, not depending on $R$, such that $\limsup_{\beta \to \infty} \frac{1}{\beta} \log \mathbb{P}_{N,\beta}^{(\otimes)}(\mu_\beta^{(i)} \notin K_R) \le -R + C$ for all $R > 0$, and this implies the exponential tightness.

The Gâteau-differentiability of $\Lambda_N^{(\otimes)}$ is proven as follows. The proof of Lemma 1.3 shows that the infimum in the formula on the right-hand side of (2.2) is attained. We abbreviate $\Lambda = \Lambda_N^{(K)}$. Fix $\Phi \equiv f = (f_1, \ldots, f_N) \in \mathcal{C}_b(\mathbb{R}^d)^N$ and some $g \in \mathcal{C}_b(\mathbb{R}^d)^N$. We want to show the existence of the limit $\lim_{t \to 0} \frac{1}{t}[\Lambda(f + tg) - \Lambda(f)]$. With $(\mu_1^{(t)}, \ldots, \mu_N^{(t)})$ a minimizer for the formula on the right-hand side of (2.2) for $f$ replaced by $f + tg$, we obtain, by replacing the minimizer in the formula for $f$ by $(\mu_1^{(t)}, \ldots, \mu_N^{(t)})$,

$$(2.33) \qquad \frac{1}{t}[\Lambda(f + tg) - \Lambda(f)] \geq \sum_{i=1}^{N} \langle g_i, \mu_i^{(t)} \rangle.$$

Since the family $(\mu_1^{(t)}, \ldots, \mu_N^{(t)})_{t>0}$ is easily seen to be convergent weakly toward the minimizer $(\mu_1, \ldots, \mu_N)$ for the formula for $\Lambda(f)$, it is clear that the right-hand side of (2.33) converges toward $\sum_{i=1}^{N} \langle g_i, \mu_i \rangle$. Analogously, one shows the complementary bound. This implies the Gâteau-differentiability of $\Lambda_N^{(K)}$ with

$$(2.34) \qquad \frac{\partial}{\partial g} \Lambda_N^{(\otimes)}(f) = \sum_{i=1}^{N} \langle g_i, \mu_i \rangle.$$

Now [11], Corollary 4.5.27, implies Theorem 1.7(ii). The statement in Theorem 1.7(iii) is a standard corollary. This ends the proof of Theorem 1.7.

## 3. Large deviations for the Dirac-interaction model: Proof of Theorem 1.12.
We follow the same strategy as in the proof of Theorem 1.7 in Section 2. Hence, the main step is the proof of the following.

PROPOSITION 3.1 (Asymptotic for the cumulant generating function). *For any* $f_1, \ldots, f_N \in \mathcal{C}_{\mathrm{b}}(\mathbb{R}^d)$,

$$
\begin{aligned}
(3.1) \qquad \lim_{\beta \to \infty} \beta^{-p/(2+p)} &\log \mathbb{E}[e^{-1/\beta \sum_{i=1}^{N} \int_0^\beta W(S_t^{(i)}) \, dt} \\
&\qquad \times e^{-\lambda \beta^{(d+p)/(2+p)} \sum_{1 \leq i < j \leq N} \alpha_\beta^{(i,j)}} e^{\beta^{p/(2+p)} \langle f, \mu_\beta^\otimes \rangle}] \\
&= -N \chi_N^{(\delta_0)}(\lambda, f),
\end{aligned}
$$

*where*

$$
\begin{aligned}
(3.2) \qquad N \chi_N^{(\delta_0)}(\lambda, f) = \inf_{h_1, \ldots, h_N \in H^1(\mathbb{R}^d) \, : \, \|h_i\|_2 = 1 \, \forall i} \Bigg[ &\sum_{i=1}^{N} I_1(h_i^2) + \langle \mathfrak{W} - f, (h^2)^\otimes \rangle \\
&+ \lambda \sum_{1 \leq i < j \leq N} \langle h_i^2, h_j^2 \rangle \Bigg],
\end{aligned}
$$

*where we wrote* $(h^2)^\otimes = h_1^2 \otimes \cdots \otimes h_N^2$ *and* $f = f_1 \oplus \cdots \oplus f_N$.

Theorem 1.12 follows from Proposition 3.1 in the same way as Theorem 1.7 follows from Proposition 2.1 in Section 2.4; we omit the details. Hence, it remains to prove Proposition 3.1, which we do in the next two sections. Let us first remark that

$$\frac{1}{\beta} \int_0^\beta W(S_t^{(i)}) \, dt = \beta^{p/(2+p)} \langle W, \mu_\beta^{(i)} \rangle \quad \text{and}$$

$$\beta^{(d+p)/(2+p)} \alpha_\beta^{(i,j)} = \beta^{p/(2+p)} \left\langle \frac{d\mu_\beta^{(i)}}{dx}, \frac{d\mu_\beta^{(j)}}{dx} \right\rangle,$$

as is derived by an elementary calculation. According to Lemma 1.16, $(\mu_\beta^{(i)})_{\beta>0}$ satisfies large deviations principle bounds with scale $\beta^{p/(2+p)}$. Hence, the problem left to be solved is to circumvent the missing boundedness and continuity of the functionals $\mu \mapsto \langle W, \mu \rangle$ and $(\mu_1, \mu_2) \mapsto \langle \frac{d\mu_1}{dx}, \frac{d\mu_2}{dx} \rangle$.

3.1. *Proof of the upper bound in* (3.1). We proceed as in Section 2.1. However, an additional smoothing argument will be necessary.

Let auxiliary parameters $\eta > 0$, $M > 0$ and $R > 0$ be given, recall $\Lambda_R = [-R, R]^d$, and consider the $\Lambda_R$-periodizations $L_{\beta,R}^{(i)}$ of the densities $L_\beta^{(i)}$ defined in (1.37). As in Section 2.1, we distinguish whether or not $\mu_\beta^{(i)}(\Lambda_R) > 1 - \eta$. Furthermore, we estimate

$$\langle W, \mu_\beta^{(i)} \rangle \geq \langle W_M, \mu_\beta^{(i)} \rangle \quad \text{and} \quad \left\langle \frac{d\mu_\beta^{(i)}}{dx}, \frac{d\mu_\beta^{(j)}}{dx} \right\rangle \geq \langle L_\beta^{(i)}, L_\beta^{(j)} \wedge M \rangle,$$

where we recall that $W_M = W \wedge M$ is the cut-off version of the trap potential. We intend to replace $L_\beta^{(i)}$ by $L_{\beta,R}^{(i)}$. The replacement errors for the terms involving $W$ and $f$ have been estimated in Section 2.1. The one for the interaction term is estimated as follows,

$$\langle L_\beta^{(i)}, L_\beta^{(j)} \wedge M \rangle - \langle L_{\beta,R}^{(i)}, L_{\beta,R}^{(j)} \wedge M \rangle$$

$$\geq \int_{\mathbb{R}^d} L_\beta^{(i)} (L_\beta^{(j)} \wedge M)(x) \, dx - \int_{\Lambda_R} L_{\beta,R}^{(i)}(x) (L_{\beta,R}^{(j)} \wedge M)(x) \, dx$$

$$= \sum_{k \in \mathbb{Z}^d \setminus \{0\}} \int_{\Lambda_R} (L_\beta^{(i)}(x + 2Rk)(L_\beta^{(j)} \wedge M)(x + 2Rk)$$

$$\qquad\qquad\qquad - L_{\beta,R}^{(i)}(x)(L_{\beta,R}^{(j)} \wedge M)(x)) \, d(x + 2Rk)$$

$$\geq -M\mu_\beta^{(i)}(\Lambda_R^c) \geq -\eta M.$$

Hence, as in Section 2.1, we obtain the bound

$$\mathbb{E}[e^{-1/\beta \sum_{i=1}^N \int_0^\beta W(S_t^{(i)}) \, dt} e^{-\lambda \beta^{(d+p)/(2+p)} \sum_{1 \leq i < j \leq N} \alpha_\beta^{(i,j)}} e^{\beta^{p/(2+p)} \langle f, \mu_\beta^\otimes \rangle}]$$

$$(3.3) \quad \leq N e^{-\beta^{p/(2+p)} (\eta \inf_{\Lambda_R^c} W - C_f)}$$

$$\qquad + e^{\beta^{p/(2+p)} \eta [NM + N^2 M + C_f]}$$

$$\qquad\qquad \times \mathbb{E}[e^{-\beta^{p/(2+p)} [\langle \mathfrak{W}_M - f, \mu_{\beta,R}^\otimes \rangle + \lambda \sum_{i<j} \langle L_{\beta,R}^{(i)}, L_{\beta,R}^{(j)} \wedge M \rangle]}].$$

The last functional in the exponent on the right-hand side is bounded, but not continuous. Hence, we need a smoothing argument. For this purpose, let $\kappa_\delta$ denote the Gaussian density in $\mathbb{R}^d$ with variance $\delta > 0$, and denote convolution by $*$. Fix some small $\varepsilon > 0$. According to [15], Lemma 3.7, specialized to our situation,

$$(3.4) \quad \lim_{\delta \downarrow 0} \limsup_{\beta \to \infty} \beta^{-p/(2+p)} \log \sup_{\|g\|_\infty \le M} \mathbb{P}(|\langle L_\beta^{(i)}, g - g * \kappa_\delta \rangle| > \varepsilon) = -\infty.$$

Actually, in [15] only the discrete-time case is handled, but the continuous-time case is similar. This means that we can estimate the last expectation on the right-hand side of (3.3) from above against

$$(3.5) \quad \begin{aligned} & e^{-\beta^{p/(2+p)} \lambda C(M,\varepsilon,\delta)} \\ & \quad + e^{-\beta^{p/(2+p)} \lambda \varepsilon} \mathbb{E}\big[ e^{-\beta^{p/(2+p)} [\langle \mathfrak{W}_M - f, \mu_{\beta,R}^\otimes \rangle - \lambda \sum_{i<j} \langle L_{\beta,R}^{(i)}, (L_{\beta,R}^{(j)} \wedge M) * \kappa_\delta \rangle ]} \big], \end{aligned}$$

where $C(M,\varepsilon,\delta)$ is a constant that satisfies $\lim_{\delta \downarrow 0} C(M,\varepsilon,\delta) = \infty$ for any $M, \varepsilon > 0$. The functional

$$(\mu_1, \mu_2) \mapsto \left\langle \frac{d\mu_1}{dx}, \frac{d\mu_2}{dx} \wedge M * \kappa_\delta \right\rangle = \left\langle \frac{d\mu_1}{dx} * \kappa_{\delta/2}, \frac{d\mu_2}{dx} \wedge M * \kappa_{\delta/2} \right\rangle$$

is bounded and continuous in the weak topology on the set of probability densities on $\Lambda_R$. Hence, we can apply Varadhan's integral lemma and the large deviation principle in Lemma 1.16(i) to the expectation in (3.5). This gives that the limit superior on the left-hand side of (3.1) is not bigger than

$$(3.6) \quad \begin{aligned} - \min \Big\{ & \eta \inf_{\Lambda_R^c} W - C_f, -\eta[NM + N^2 M + C_f] + \lambda C(M,\varepsilon,\delta), \\ & -\eta N[NM + N^2 M + C_f] - \lambda \varepsilon - N \chi_N^{(\delta_0)}(R,M,\delta,\lambda,f) \Big\}, \end{aligned}$$

where

$$\begin{aligned} \chi_N^{(\delta_0)} & (R,M,\delta,\lambda,f) \\ & = \inf \Big\{ \sum_{i=1}^N I_R^{(\mathrm{per})}(\mu_i) + \langle \mathfrak{W}_M - f, \mu_1 \otimes \cdots \otimes \mu_N \rangle \\ & \qquad + \sum_{i<j} \left\langle \frac{d\mu_i}{dx}, \left( \frac{d\mu_j}{dx} \wedge M \right) * \kappa_\delta \right\rangle : \\ & \qquad\qquad\qquad \mu_1, \ldots, \mu_N \in \mathcal{M}_1(\Lambda_R) \Big\}. \end{aligned}$$

Now we let $\delta \downarrow 0$, $\varepsilon \downarrow 0$, $R \to \infty$, $\eta \downarrow 0$ and finally $M \to \infty$. It is elementary to derive that $\liminf_{\delta \downarrow 0} \chi_N^{(\delta_0)}(R,M,\delta,\lambda,f) \ge \chi_N^{(\delta_0)}(R,M,0,\lambda,f)$, where we interpret $\mu * \kappa_0$ as $\mu$. Furthermore, similarly to the proof of (2.14), one shows that

$$\liminf_{M \to \infty} \liminf_{R \to \infty} \chi_N^{(\delta_0)}(R, M, 0, \lambda, f) \geq \chi_N^{(\delta_0)}(\lambda, f).$$

This ends the proof of the upper bound in (3.1).

3.2. *Proof of the lower bound in* (3.1). For the proof of the lower bound, we follow a discrete variant of the strategy used in Section 2.3, that, we use $N$ eigenvalue expansions for the $N$ expectations over the $N$ random walks separately. First we consider the expectation with respect to the $N$th walk, almost surely with respect to the other $N-1$ walks. We only treat the terms depending on the $N$th motion and introduce the potential

$$q_N(z) = \xi_\beta^{-2} \left[ W\left(\frac{z}{\xi_\beta}\right) + \lambda \sum_{j<N} L_\beta^{(j)}\left(\frac{z}{\xi_\beta}\right) - \overline{f}_N\left(\frac{z}{\xi_\beta}\right) \right],$$

where $\xi_\beta = \beta^{1/(2+p)}$, and $\overline{f}_N(x) = \xi_\beta^d \int_{x+[0,\xi_\beta^{-1}]^d} f_N(y)\,dy$. Denote by $\mathbb{E}^{(i)}$ the expectation with respect to the $i$th random walk $S^{(i)}$ and by $\nu_i$ its initial distribution. Choose $R > 0$ so large that $\nu_i(B_{(R-1)\xi_\beta}) > 0$. Then we have

$$(3.7) \quad \begin{aligned} &\mathbb{E}^{(N)}\left[e^{-1/\beta \int_0^\beta W(S_t^{(N)})\,dt - \lambda\beta^{(d+p)/(2+p)} \sum_{j<N} \alpha_\beta^{j,N} + \beta^{p/(2+p)} \langle f_N, L_\beta^{(N)} \rangle}\right] \\ &\geq \sum_{z \in B_{R\xi_\beta}} \nu_N(z) \mathbb{E}_z^{(N)}\left[e^{-\int_0^\beta q_N(S_t^{(N)})\,dt} \mathbb{1}\{\mathrm{supp}(L_\beta^{(N)}) \subset B_{R\xi_\beta}\} \mathbb{1}\{S_\beta^{(N)} = z\}\right] \\ &\geq \sum_{z \in B_{R\xi_\beta}} \nu_N(z) e_N(z)^2 e^{\beta\lambda_N}, \end{aligned}$$

where $\lambda_N$ is the principal eigenvalue of $\Delta - q_N$ in $B_{R\xi_\beta}$ with zero boundary condition, and $e_N$ is the corresponding positive $\ell^2$-normalized eigenfunction. Pick some $L^2$-normalized function $h_N \in \mathcal{C}^2(\mathbb{R}^d)$ satisfying $\mathrm{supp}(h_N) \subset \Lambda_R$ and pick $v_N(z) = (1 + o(1))\xi_\beta^{-d/2} h_N(z\xi_\beta^{-1})$, where $o(1)$ is chosen such that $v_N$ is $\ell^2$-normalized. Then we have, using the Rayleigh–Ritz principle for the principal eigenvalue, and noting that $\overline{f}_N \to f_N$ as $\beta \to \infty$ uniformly,

$$\begin{aligned} \lambda_N &\geq \left[ -\|\nabla v_N\|_2^2 - \xi_\beta^{-2} \left\langle W(\cdot \xi_\beta^{-1}) - \overline{f}_N(\cdot \xi_\beta^{-1}) + \lambda \sum_{j<N} L_\beta^{(j)}(\cdot \xi_\beta^{-1}), v_N^2(\cdot) \right\rangle \right] \\ &\quad \times (1 + o(1)) \\ &= -\xi_\beta^{-2} \left[ \|\nabla h_N\|_2^2 + \langle W - \overline{f}_N, h_N^2 \rangle + \lambda \sum_{j<N} \langle L_\beta^{(j)}, h_N^2 \rangle \right] + o(\xi_\beta^{-2}). \end{aligned}$$

We estimate the term $\sum_z \nu_N(z) e_N(z)^2$ on the right-hand side of (3.7). Note that $e_N$ is also an eigenfunction for the transition densities of the random walk in $B_{R\xi_\beta}$ with potential $-q_N - \lambda_N$, that is,

$$(3.8) \quad e_N(z) = \mathbb{E}_z[e^{-\int_0^1 (q_N(S_t) + \lambda_N)\,dt} \mathbb{1}\{S_t \in B_{R\xi_\beta} \,\forall t \in [0,1]\} e_N(S_1)],$$
$$z \in B_{R\xi_\beta},$$

where $\mathbb{E}_z$ is the expectation with respect to an independent copy $(S_t)_{t\in[0,\infty)}$ of, say, $(S_t^{(1)})_{t\in[0,\infty)}$. We can estimate $\lambda_N \leq \xi_\beta^{-2} \|f_N\|_\infty$ and, since $L_\beta^{(j)}(x) \leq \xi_\beta^d$ for any $x \in \mathbb{R}^d$ and any $j \in \{1,\ldots,N\}$,

$$q_N(S_t) \leq \xi_\beta^{-2}\left[R^p + \lambda \sum_{j<N} L_\beta^{(j)}(S_t) + \|f_N\|_\infty\right]$$

$$\leq [R^p + \|f_N\|_\infty]\beta^{-2/(2+p)} + \lambda(N-1)\beta^{(d-2)/(2+p)} \leq C\beta^{(d-2)/(2+p)},$$

for some $C > 0$, depending only on $R$, $\lambda$, $N$, $p$ and $\|f_N\|_\infty$. Using this in (3.8), we find, for all large $\beta$ (changing the value of $C$ if necessary),

$$e_N(z) \geq e^{-C\beta^{(d-2)/(2+p)}} \sum_{\widetilde{z} \in B_{R\xi_\beta}} \mathbb{P}_z(S_t \in B_{R\xi_\beta} \,\forall t \in [0,1], S_1 = \widetilde{z}) e_N(\widetilde{z}),$$
$$z \in B_{R\xi_\beta}.$$

Recall that $p > d - 2$, hence, $\beta^{(d-2)/(2+p)} = o(\beta^{p/(2+p)})$. It is clear that $\mathbb{P}_z(S_t \in B_{R\xi_\beta} \,\forall t \in [0,1], S_1 = \widetilde{z}) \geq e^{-o(\beta^{p/(2+p)})}$, uniformly in $z, \widetilde{z} \in B_{R\xi_\beta}$. Since $e_N$ is $\ell^2$-normalized, we can estimate $\sum_{\widetilde{z} \in B_{R\xi_\beta}} e_N(\widetilde{z}) \geq \sum_{\widetilde{z} \in B_{R\xi_\beta}} e_N(\widetilde{z})^2 = 1$. Pick some $z_N \in B_{R\xi_\beta}$ such that $\nu_N(z_N) > 0$, then we have from the preceding, for any large $\beta > 0$,

$$(3.9) \quad \begin{aligned} &\mathbb{E}^{(N)}[e^{-1/\beta \int_0^\beta W(S_t^{(N)})\,dt - \lambda\beta^{(d+p)/(2+p)}\sum_{j<N}\alpha_\beta^{j;N} + \beta^{p/(2+p)}\langle f_N, L_\beta^{(N)}\rangle}] \\ &\geq \exp\left\{-\beta^{p/(2+p)}\left[\|\nabla h_N\|_2^2 + \langle W - f_N, h_N^2\rangle + \sum_{j<N}\langle L_\beta^{(j)}, h_N^2\rangle\right] \right. \\ &\left. \qquad\qquad + o(\beta^{p/(2+p)})\right\}, \end{aligned}$$

and the term $o(\beta^{p/(2+p)})$ does not depend on the $N - 1$ other random walks.

In the same way, we treat the expectations with respect to the other $N-1$ random walks. For doing this, introduce recursively, for $i = N-1, N-2, \ldots, 1$, the potentials

$$q_i(z) = \xi_\beta^{-2}\left[W\left(\frac{z}{\xi_\beta}\right) + \lambda \sum_{j<i} L_\beta^{(j)}\left(\frac{z}{\xi_\beta}\right) - \sum_{j>i} h_j^2\left(\frac{z}{\xi_\beta}\right) - \overline{f}_i\left(\frac{z}{\xi_\beta}\right)\right],$$

where $\overline{f}_i(x) = \xi_\beta^d \int_{x+[0,\xi_\beta^{-1}]^d} f_i(y)\,dy$, and pick $L^2$-normalized functions $h_{N-1}$, $h_{N-2}, \ldots, h_1 \in \mathcal{C}^2(\mathbb{R}^d)$ with supports in $\Lambda_R$. In the same way as (3.9), one derives

$$\mathbb{E}[e^{-1/\beta \sum_{i=1}^N \int_0^\beta W(S_t^{(i)})\,dt} e^{-\lambda \beta^{(d+p)/(2+p)} \sum_{1 \le i < j \le N} \alpha_\beta^{(i,j)}} e^{\beta^{p/(2+p)} \langle f, \mu_\beta^\otimes \rangle}]$$

$$\ge e^{o(\beta^{p/(2+p)})} \exp\left\{-\beta^{p/(2+p)} \left[\sum_{i=1}^N [\|\nabla h_i\|_2^2 + \langle W - f_i, h_i^2 \rangle] \right.\right.$$
$$\left.\left. - \lambda \sum_{1 \le i < j \le N} \langle h_i^2, h_j^2 \rangle \right]\right\}.$$

Picking the $L^2$-normalized functions $h_1, \ldots, h_N$ with support in $\Lambda_R$ optimally and letting $R \to \infty$, we arrive the lower bound in (3.1), noting the substitution $h_i^2(x)\,dx = \mu_i(dx)$ in (3.2).

## 4. Analysis of the ground product states.
We prove Lemma 1.3 in Sections 4.1–4.4 and Theorem 1.14 in Section 4.5. The proofs combine methods from variational analysis and from probabilistic potential theory.

### 4.1. *Existence of minimizers in* (1.11).
Let $(h_{1,k}, \ldots, h_{N,k})_{k \in \mathbb{N}}$ be a sequence of approximate minimizers for the right-hand side of (1.11), that is, $h_{i,k} \in H^1(\mathbb{R}^d)$ with $\|h_{i,k}\|_2 = 1$ for all $i = 1, \ldots, N$ and $k \in \mathbb{N}$, and $N\chi_N^{(\otimes)} = \lim_{k \to \infty} [\sum_{i=1}^N \|\nabla h_{i,k}\|_2^2 + \langle \mathfrak{W} + \mathfrak{v}, h_k^2 \rangle]$, where we abbreviated $h_k = h_{1,k} \otimes \cdots \otimes h_{N,k}$. Fix $i \in \{1, \ldots, N\}$. Our first goal is to establish the convergence of $(h_{i,k})_{k \in \mathbb{N}}$ along suitable subsequences.

For any $R > 0$ and any $k \in \mathbb{N}$, we have

$$\begin{aligned}
\int_{\Lambda_R^c} h_{i,k}^2(x)\,dx &\le \frac{\int_{\mathbb{R}^d} h_{i,k}^2(x) W(x)\,dx}{\inf_{\Lambda_R^c} W} \\
&\le \frac{\langle \mathfrak{W}, h_k^2 \rangle}{\inf_{\Lambda_R^c} W} \\
&\le \frac{\sum_{i=1}^N \|\nabla h_{i,k}\|_2^2 + \langle \mathfrak{W} + \mathfrak{v}, h_k^2 \rangle - ((N(N-1))/2)\inf v}{\inf_{\Lambda_R^c} W}.
\end{aligned}$$

(4.1)

By our assumption in (1.1), $\sup_{k \in \mathbb{N}} \int_{\Lambda_R^c} h_{i,k}^2(x)\,dx$ vanishes as $R \to \infty$. Hence, the family of probability measures $(h_{i,k}^2(x)\,dx)_{k \in \mathbb{N}}$ is tight, and from Prohorov's theorem, we conclude the existence of a probability measure $\mu_i \in \mathcal{M}_1(\mathbb{R}^d)$ such that $(h_{i,k}^2(x)\,dx)_{k \in \mathbb{N}}$ converges weakly toward $\mu_i$ along a subsequence, which we again denote $(k)_{k \in \mathbb{N}}$. Hence, for any bounded continuous function $f : \mathbb{R}^d \to \mathbb{R}$, we have $\lim_{k \to \infty} \langle f, h_{i,k}^2 \rangle = \langle f, \mu_i \rangle$.

Furthermore, because $\|h_{i,k}\|_2 = 1$ for any $k \in \mathbb{N}$, according to Banach–Alaoglu's theorem, there is a $h_i \in L^2(\mathbb{R}^d)$ such that $(h_{i,k})_{k \in \mathbb{N}}$ converges weakly in $L^2$ along a subsequence [which we again denote $(k)_{k\in\mathbb{N}}$] toward $h_i$. Since, by approximative minimality and because $W$ and $v$ are bounded from below, also the sequence of energies $(\|\nabla h_{i,k}\|_2)_{k\in\mathbb{N}}$ is bounded, the Banach–Alaoglu theorem implies that also the sequence of gradients converge weakly in $L^2$ along a subsequence [which we again denote $(k)_{k\in\mathbb{N}}$] toward some $g_i$. In particular, $h_i \in H^1(\mathbb{R}^d)$ and $\nabla h_i = g_i$, which follows via [21], Theorem 8.6, from the completeness of the Sobolov space. By [21], Theorem 8.7, we may assume that the convergence $h_{i,k} \to h_i$ is pointwise almost everywhere.

We argue now that $h_i^2$ is a density of $\mu_i$. According to [21], Theorem 8.6, the convergence of $(h_{i,k})_k$ toward $h_i$ is also strong in $L^2$ on every compact set. In particular, the integrals of $h_{i,k}^2$ against any continuous function with compact support converges to that of $h_i^2$. By weak convergence of $h_{i,k}^2(x)\,dx$ toward $\mu_i$, they converge also toward the integrals against $\mu_i$. Hence, the integrals of any bounded, compactly supported function against $h_i^2$ and against $\mu_i$ coincide, which implies that $\mu_i(dx) = h_i^2(x)\,dx$.

Now we prove that the vector $(\mu_1,\ldots,\mu_N)$ is a minimizer on the right-hand side of (1.11). To do this, we verify suitable lower-semicontinuities of the functional on the right-hand side of (1.11). It is well known that the energy functional, $h \mapsto \|\nabla h\|_2$, is weakly lower-semicontinuous in $L^2$. Furthermore, the map $\mu \mapsto \langle W, \mu \rangle$ is obviously lower-semicontinuous in the weak topology of probability measures, since $W$ is nonnegative and continuous. By $\mathfrak{v}_+$ and $\mathfrak{v}_-$ we denote the positive and the negative part of $\mathfrak{v}$, respectively. Fatou's lemma implies that $\liminf_{k\to\infty}\langle \mathfrak{v}_+, h_k^2 \rangle \geq \langle \mathfrak{v}_+, h^2 \rangle$, where we recall that $h_k = h_{1,k} \otimes \cdots \otimes h_{N,k}$. Furthermore, note that

$$|\langle \mathfrak{v}_+, h_k^2 \rangle - \langle \mathfrak{v}_+, h^2 \rangle|$$
$$\leq -\inf v \sum_{i<j} \int_{\mathbb{R}^d} \int_{\mathbb{R}^d} |h_{i,k}^2(x) h_{j,k}^2(x+\widetilde{x}) - h_i^2(x) h_j^2(x+\widetilde{x})|\,dx\,d\widetilde{x}.$$

Using (4.1) and the strong $L^2$ convergence on compacts, one easily derives that the right-hand side vanishes as $k \to \infty$. Summarizing,

$$\sum_{i=1}^N \|\nabla h_i\|_2^2 + \langle \mathfrak{W} + \mathfrak{v}, h^2 \rangle \leq \liminf_{k\to\infty} \left[ \sum_{i=1}^N \|\nabla h_{i,k}\|_2^2 + \langle \mathfrak{W} + \mathfrak{v}, h_k^2 \rangle \right] = N\chi_N^{(\otimes)},$$

where we wrote $h = h_1 \otimes \cdots \otimes h_N$ and (in an abuse of notation) $h_k = h_{1,k} \otimes \cdots \otimes h_{N,k}$. Hence, the tuple $(h_1,\ldots,h_N)$ is a minimizer of the right-hand side of (1.11).

### 4.2. *Positivity and regularity of minimizers, soft-core case.*

First we consider the case of a soft-core potential $v$. Let $(h_1,\ldots,h_N) \in H^1(\mathbb{R}^d)^N$ be a minimizer of the variational formula in (1.9). Fix $i \in \{1,\ldots,N\}$. We prove that $h_i$ is strictly positive everywhere in $\mathbb{R}^d$. We certainly may assume that $h_i \geq 0$. Note that $h_i$ is a minimizer in the problem

$$(4.2) \quad \lambda_i = \inf_{\widetilde{h}_i \in H^1(\mathbb{R}^d), \|\widetilde{h}_i\|_2 = 1} \left[ \tfrac{1}{2} \|\nabla \widetilde{h}_i\|_2^2 + \langle W, \widetilde{h}_i^2 \rangle + \left\langle \widetilde{h}_i^2, \sum_{j \neq i} V h_j^2 \right\rangle \right].$$

First we show that $V h_j^2$ is locally integrable for any $j \in \{1, \dots, N\}$ with $j \neq i$. Indeed, for any measurable set $B \subset \mathbb{R}^d$ of finite measure, we have [denoting the $\varepsilon$-ball around $x$ by $B_\varepsilon(x)$]

$$
\begin{aligned}
(4.3) \quad \int_B V h_j^2(x) \, dx &= \int_{\mathbb{R}^d} dy \, v(|y|) \int_B dx \, h_j^2(y + x) \\
&\leq \int_{B_\varepsilon(0)} dy \, v(|y|) \|h_j\|_2^2 + \sup_{[\varepsilon, \infty)} v \int_B dx \int_{B_\varepsilon^c(0)} dy \, h_j^2(x + y) \\
&\leq \int_{B_\varepsilon(0)} dy \, v(|y|) + |B| \sup_{[\varepsilon, \infty)} v \\
&< \infty.
\end{aligned}
$$

Since $V h_j^2$ is also bounded from below, the interaction potential of the problem in (4.2) is locally integrable. The same holds for the trap part, thus, the potential $W + \sum_{j \neq i} V h_j^2 - \lambda_i$ in (4.2) is locally integrable. By [21], Theorem 11.8, $h_i$ satisfies the Euler–Lagrange equation

$$(4.4) \quad \Delta h_i = U_i h_i \quad \text{in } \mathcal{D}'(\mathbb{R}^d), \qquad \text{where } U_i = W + \sum_{j \neq i} V h_j^2 - \lambda_i.$$

(Note here that the potential $U_i$ is a priori not locally bounded, but the first part of the proof of [21], Theorem 11.8, does not use this.) So, we get by [21], Theorem 11.7, that

$$(4.5) \quad h_i \in L^q_{\text{loc}}(\mathbb{R}^d) \qquad \text{for all } q < \begin{cases} \infty, & \text{in } d = 2, \\ 3, & \text{in } d = 3. \end{cases}$$

Now proceed for $d = 2$ and $d = 3$ separately. In $d = 2$ we deduce that $V h_j^2 \in L^p_{\text{loc}}(\mathbb{R}^2)$ for any $p \in [1, \infty)$. Indeed, fix any bounded measurable set $B \subset \mathbb{R}^d$, pick $\eta > 0$ so small that $v(|y|) \geq 0$ for all $y \in B_\eta(0)$ and abbreviate $K_\eta = \int_{B_\eta(0)} dy \, v(|y|)$. Then, using Jensen's inequality and (4.5), we see that

$$
\begin{aligned}
(4.6) \quad \int_B |V h_j^2(x)|^p \, dx &\leq K_\eta^p \int_B dx \left| \int_{B_\eta(x)} dy \, \frac{v(|x - y|)}{K_\eta} h_j^2(y) \right|^p + |B| \sup_{(\eta, \infty)} |v|^p \\
&\leq K_\eta^p \int_B dx \int_{B_\eta(x)} dy \, \frac{v(|x - y|)}{K_\eta} h_j(y)^{2p} + |B| \sup_{(\eta, \infty)} |v|^p \\
&= K_\eta^{p-1} \int_{B_\eta(0)} dy \, v(|y|) \int_B dx \, h_j^{2p}(y + x) + |B| \sup_{(\eta, \infty)} |v|^p \\
&< \infty.
\end{aligned}
$$

Hence, the potential $U_i$ is locally $p$-integrable for any $p < \infty$. Using again [21], Theorem 11.7, we obtain that $h_i$ is continuously differentiable, and all first partial derivatives of $h_i$ are $\alpha$-Hölder continuous for any $\alpha < 1$. In particular, $U_i$ is also locally bounded, and from [21], Theorem 9.10, we have $h_i > 0$ everywhere in $\mathbb{R}^2$. For this, we apply [21], Theorem 9.10, to any open, bounded and connected $A \subset \mathbb{R}^2$, that satisfies $\int_A h_i(x)\, dx > 0$.

Now we turn to the case $d = 3$. Recall that we assume that

$$\int_{B_1(0)} |v(|y|)|^{1+\delta}\, dy < \infty \qquad \text{for some } \delta > 0.$$

We show that $|Vh_j^2|^p$ is locally integrable for some $p > \frac{3}{2}$. Recall that $h_j^2 \in L^{p'}_{\mathrm{loc}}(\mathbb{R}^3)$ for any $p' < \frac{3}{2}$. Pick a bounded open ball $\widetilde{B} \subset \mathbb{R}^3$ such that $\int_{\widetilde{B}} h_j^{2p'}(x)\, dx < \infty$. Pick some open ball $B$ whose closure is contained in $\widetilde{B}$. We shall show that $\int_B |Vh_j^2(x)|^p\, dx < \infty$ for some $p > \frac{3}{2}$. Let now $\eta > 0$ be so small that $B_\eta(x) \subset \widetilde{B}$ for any $x \in B$ and that $v(|y|) \geq 0$ for all $y \in B_\eta(0)$. We estimate $|Vh_j^2(x)| \leq \sup_{(\eta,\infty)} |v| + \int_{B_\eta(0)} dy\, v(|y|)h_j^2(x+y)$ for $x \in B$. Hence, we have to show that

$$(4.7) \qquad \int_B dx \left| \int_{B_\eta(0)} dy\, v(|y|)h_j^2(x+y) \right|^p < \infty \qquad \text{for some } p > \frac{3}{2}.$$

We shall do that for any $p > 3/2$ and $p' < 3/2$ satisfying $p/p' < 1 + \delta$. Recall the abbreviation $K_\eta = \int_{B_\eta(0)} dy\, v(|y|)$. Jensen's inequality gives

$$
\begin{aligned}
\int_B dx &\left| \int_{B_\eta(0)} dy\, v(|y|)h_j^2(x+y) \right|^p \\
&\leq K_\eta^p \int_B dx \left| \int_{B_\eta(0)} dy\, \frac{v(|y|)}{K_\eta} h_j^{2p'}(x+y) \right|^{p/p'} \\
&= K_\eta^{p-p/p'} \int_B dx \left| \int_{B_\eta(0)} dy \left( \frac{h_j^{2p'}(x+y)}{\int_{B_\eta(x)} dz\, h_j^{2p'}(z)} \right) v(|y|) \right|^{p/p'} \\
&\qquad \times \left( \int_{B_\eta(x)} dz\, h^{2p'}(z) \right)^{p/p'} \\
&\leq K_\eta^{p-p/p'} \int_B dx \int_{B_\eta(0)} dy\, h_j^{2p'}(x+y) v(|y|)^{p/p'} \\
&\qquad \times \left( \int_{B_\eta(x)} dz\, h_j^{2p'}(z) \right)^{p/p'-1} \\
&< \infty,
\end{aligned}
$$

(4.8)

where we used that $p < p'(1+\delta)$ in the last step. Hence, $U_i \in L^p_{\mathrm{loc}}(\mathbb{R}^3)$ for some $p > \frac{3}{2}$. Now [21], Theorem 11.7, implies that $h_i$ is continuous and that $Vh_i^2$ is locally bounded. Therefore, $U_i$ is locally bounded. Finally, we get from [21], Theorem 11.7, that $h_i$ is continuously differentiable and that the partial derivatives are $\alpha$-Hölder continuous for any $\alpha < 1$; whereas [21], Theorem 9.10, gives that $h_i$ is positive everywhere in $\mathbb{R}^3$.

4.3. *Positivity and regularity of minimizers, hard-core case.* Now let $v$ be a hard-core potential, and recall the parameter $a \in [0, \infty)$ from (1.2). Let $(h_1, \ldots, h_N)$ be a minimizer of the formula in (1.11). Fix $i \in \{1, \ldots, N\}$ and put

$$\Omega_i = \overline{B_a\left(\bigcup_{j \neq i} \mathrm{supp}(h_j)\right)}^{\mathrm{c}},$$

where $B_a(A)$ is the $a$-neighborhood of a set $A$ if $a > 0$, and $B_0(A) = A$. Then it is easy to see that, for any $j \neq i$, the potential $Vh_j^2$ is bounded on any compact subset of $\Omega_i$. It is clear that $\mathrm{supp}(h_i) \subset \Omega_i$. Hence, $h_i$ is a minimizer in the problem

$$(4.9) \qquad \inf_{\widetilde{h}_i \in H^1(\mathbb{R}^d), \, \mathrm{supp}(\widetilde{h}_i) \subset \Omega_i, \|\widetilde{h}_i\|_2 = 1} \left[ \|\nabla \widetilde{h}_i\|_2^2 + \langle W, \widetilde{h}_i^2 \rangle + \left\langle \widetilde{h}_i^2, \sum_{j \neq i} Vh_j^2 \right\rangle \right].$$

Again by [21], Theorem 11.8, $h_i$ satisfies the Euler–Lagrange equation in (4.4) in $\mathcal{D}'(\Omega_i)$. By [21], Theorem 11.7, $h_i$ is continously differentiable in $\Omega_i$, and its first partial derivatives are $\alpha$-Hölder continuous for any $\alpha < 1$. By [21], Theorem 9.10, $h_i$ is strictly positive throughout any open component of $\Omega_i$ in which it is not identically equal to zero.

4.4. *Boundedness of minimizers.* Now we prove that $\|h_i\|_\infty \leq C_d(\lambda_i - (N-1)\inf v)^{d/4}$. We do this by carefully examining the proof of [3], Corollary 2.5. Recall the potential $U_i = W + \sum_{j \neq i} Vh_j^2 - \lambda_i$ from (4.4), which is locally bounded and globally bounded from below, more precisely, $-\inf U_i \leq \lambda_i - (N-1)\inf v =: u_i \in [0, \infty)$. Fix $x_0 \in \mathbb{R}^d$ and $0 < r < (4Cu_i)^{-1/2}$, where $C = \sup_{|x| < 1} \mathbb{E}_x(T_1(0))$, and $T_r(y)$ denotes the first exit time of a Brownian motion from the ball $B_r(y)$. Then we have $\mathbb{E}_{x_0}(2u_iT_r(x_0)) = 2u_iCr^2 \leq \frac{1}{2}$ and Khas'minskii's lemma ([3], Theorem 1.2) implies that

$$(4.10) \qquad \mathbb{E}_{x_0}\left( \exp\left\{ -\int_0^{T_r(x_0)} 2U_i(B_s)\, ds \right\} \right) \leq \mathbb{E}_{x_0}(e^{2u_iT_r(x_0)}) \leq 2.$$

Note that $h_i(x) = \mathbb{E}_x[e^{-\int_0^{T_r(x_0)} U_i(B(s))\, ds} h_i(B(T_r(x_0)))]$ for $x \in B_r(x_0)$ by [3], Theorem A.4.1. We apply the Cauchy–Schwarz inequality and (4.10) to get that

$$(4.11) \qquad h_i(x_0)^2 \leq 2\mathbb{E}_{x_0}[h_i^2(B(T_r(x_0)))] = 2\int_{|x-x_0|=r} \sigma(dx)h_i(x)^2,$$

where $\sigma$ denotes the normalized surface measure on $\partial B_r(x_0)$. Averaging the left- and the right-hand side of (4.11) on $r \in [0, (4Cu_i)^{-1/2}]$, we obtain

$$h_i(x_0)^2 \leq 2(4Cu_i)^{d/2} \int_{B_{(4Cu_i)^{-1/2}}(x_0)} h_i(x)^2 \, dx \leq 2(4C[\lambda_i - (N-1)\inf v])^{d/2},$$

and this implies the assertion.

This completes the proof of Lemma 1.3.

4.5. *Proof of Theorem* 1.14. Recall that we are in $d \in \{2,3\}$ and assume that $v \geq 0$ and $v(0) > 0$ and $\widetilde{\alpha}(v) < \infty$, where $\widetilde{\alpha}(v)$ is given in (1.41). Furthermore, recall that we replace $v$ in (1.9) by $v_N(\cdot) = N^{d-1}v(\cdot N)$.

First we show the upper bound, $\limsup_{N\to\infty} \chi_N^{(\otimes)} \leq \chi_{\widetilde{\alpha}(v)}^{(\mathrm{GP})}$. Let $\phi$ denote the minimizer in the Gross–Pitaevskii formula (1.40) with $\alpha = \widetilde{\alpha}(v)$. Picking $(h_1, \ldots, h_N) = (\phi, \ldots, \phi)$ in (1.11), we obtain that

$$
\begin{aligned}
(4.12) \quad \chi_N^{(\otimes)} &\leq \|\nabla\phi\|_2^2 + \langle W, \phi^2 \rangle \\
&\quad + \tfrac{1}{2} \int_{\mathbb{R}^d} N^d v(N|x|)(\phi^2 \overline{\divideontimes} \phi^2)(x) \, dx \\
&= \|\nabla\phi\|_2^2 + \langle W, \phi^2 \rangle \\
&\quad + \tfrac{1}{2} \int_{\mathbb{R}^d} v(|x|)(\phi^2 \overline{\divideontimes} \phi^2)(x/N) \, dx,
\end{aligned}
$$

where $h \overline{\divideontimes} g(x) = \int_{\mathbb{R}^d} h(y)g(x+y) \, dy$. Since $\phi^2$ is bounded, integrable and continuous, we can use Lebesgue's theorem to see that the last term converges toward $\tfrac{1}{2} \int v(|x|) \, dx (\phi^2 \overline{\divideontimes} \phi^2)(0) = 4\pi\widetilde{\alpha}(v)\|\phi\|_4^4$. Therefore, we have

$$\limsup_{N\to\infty} \chi_N^{(\otimes)} \leq \|\nabla\phi\|_2^2 + \langle W, \phi^2 \rangle + 4\pi\widetilde{\alpha}(v)\|\phi\|_4^4 = \chi_{\widetilde{\alpha}(v)}^{(\mathrm{GP})}.$$

Now we show the reversed inequality, $\liminf_{N\to\infty} \chi_N^{(\otimes)} \geq \chi_{\widetilde{\alpha}(v)}^{(\mathrm{GP})}$. Let, for any $N \in \mathbb{N}$, $(h_1^{(N)}, \ldots, h_N^{(N)})$ be any minimizer in (1.11), and define

$$\phi_N^2 = \frac{1}{N}\sum_{i=1}^N (h_i^{(N)})^2.$$

Since, for every $R > 0$,

$$(4.13) \qquad \int_{\Lambda_R^c} \phi_N^2(x) \, dx \leq \frac{\langle W, \phi_N^2 \rangle}{\inf_{\Lambda_R^c} W} \leq \frac{\chi_N^{(\otimes)}}{\inf_{\Lambda_R^c} W},$$

the sequence of probability measures $(\phi_N^2(x) \, dx)_{n\in\mathbb{N}}$ is tight [recall (1.1)]. Hence, there is a limiting probability measure $\mu$ of $\phi_N^2(x) \, dx$ in the weak sense, along some subsequence. Since $(\phi_N)_N$ is bounded in $L^2(\mathbb{R}^d)$, there is an $L^2(\mathbb{R}^d)$-weakly converging subsequence toward some $\phi \in H^1(\mathbb{R}^d)$. Since also $(\|\nabla\phi_N\|_2)_N$ is bounded [which is derived similarly as in (4.13) with the help of the convexity of the map $\phi^2 \mapsto \|\nabla\phi\|_2$], we may assume that $\nabla\phi_N \to \nabla\phi$ weakly in $L^2(\mathbb{R}^d)$. Furthermore, we may assume that the convergence $\phi_N \to \phi$ is even strong on every compact subset of $\mathbb{R}^d$ and pointwise almost everywhere. In particular, all integrals of $\phi_N$ against continuous, compactly supported functions converge to their integral against $\phi$. Hence, $\mu(dx) = \phi(x)^2\,dx$, and we have that $\|\phi\|_2 = 1$. By $L^2$-weak lower semicontinuity and convexity of the energy, we have

$$\|\nabla\phi\|_2^2 \le \liminf_{N\to\infty} \|\nabla\phi_N\|_2^2 \le \liminf_{N\to\infty} \frac{1}{N}\sum_{i=1}^{N} \|\nabla h_i^{(N)}\|_2^2$$

and

$$\langle W, \phi^2 \rangle \le \liminf_{N\to\infty} \langle W, \phi_N^2 \rangle.$$

Hence, we have

$$\liminf_{N\to\infty} \frac{1}{N} \langle h_1^{(N)} \otimes \cdots \otimes h_N^{(N)}, \mathcal{H}_N h_1^{(N)} \otimes \cdots \otimes h_N^{(N)} \rangle$$

$$\ge \frac{1}{2} \liminf_{N\to\infty} \int_{\mathbb{R}^d} v(|x|)(\phi_N^2 \,\overline{*}\, \phi_N^2)\left(\frac{x}{N}\right) dx + \|\nabla\phi\|_2^2 + \langle W, \phi^2 \rangle$$

$$- \frac{1}{2} \limsup_{N\to\infty} \frac{1}{N^2}\sum_{i=1}^{N} \int_{\mathbb{R}^d}\int_{\mathbb{R}^d} h_i^{(N)}(y)^2 h_i^{(N)}\left(y + \frac{x}{N}\right)^2 v(|x|)\,dx\,dy.$$

Hence, for proving that $\liminf_{N\to\infty} \chi_N^{(\otimes)} \ge \chi_{\overline{\alpha}(v)}^{(\mathrm{GP})}$, it remains to show that

$$\liminf_{N\to\infty} \int_{\mathbb{R}^d} v(|x|)(\phi_N^2 \,\overline{*}\, \phi_N^2)\left(\frac{x}{N}\right) dx$$

$$\ge \int_{\mathbb{R}^d} v(|x|)\,dx\|\phi\|_4^4, \tag{4.14}$$

$$\limsup_{N\to\infty} \frac{1}{N^2}\sum_{i=1}^{N} \int_{\mathbb{R}^d}\int_{\mathbb{R}^d} h_i^{(N)}(y)^2 h_i^{(N)}\left(y + \frac{x}{N}\right)^2 v(|x|)\,dx\,dy$$

$$= 0. \tag{4.15}$$

For the proof of (4.15), recall from Lemma 1.3 that $\|h_i^{(N)}\|_\infty \le C_d(\lambda_i^{(N)})^{d/4}$, where $C_d > 0$ depends only on $d$, and $\lambda_i^{(N)}$ is the $\lambda_i$ of Lemma 1.3(ii) with $v$ replaced by $v_N$. Hence,

$$\frac{1}{N^2}\sum_{i=1}^{N} \int_{\mathbb{R}^d}\int_{\mathbb{R}^d} h_i^{(N)}(y)^2 h_i^{(N)}\left(y + \frac{x}{N}\right)^2 v(|x|)\,dx\,dy$$

$$\leq C_d^2 \frac{1}{N^2} \sum_{i=1}^N (\lambda_i^{(N)})^{d/2} \int_{\mathbb{R}^d} v(|x|)\, dx.$$

Recall that $\frac{1}{N}\sum_{i=1}^N \lambda_i^{(N)} = \chi_N^{\otimes} + \frac{1}{N}\sum_{i<j}\langle (h_i^{(N)})^2, V_N(h_j^{(N)})^2\rangle \leq 2\chi_N^{\otimes}$, which is bounded from above in $N$, as we have seen in the first part of the proof. Hence, (4.15) directly follows in $d=2$, and in $d=3$ we estimate $(\lambda_i^{(N)})^{1/2} \leq CN^{1/2}$ for some $C>0$ and all $N\in\mathbb{N}$, and (4.15) also follows.

We turn to the proof of (4.14). Fatou's lemma gives that

$$\liminf_{N\to\infty} \int_{\mathbb{R}^d} v(|x|)(\phi_N^2 \,\overline{\ast}\, \phi_N^2)\left(\frac{x}{N}\right) dx$$
$$\geq \int_{\mathbb{R}^d} dx\, v(|x|) \int_{\mathbb{R}^d} dy \liminf_{N\to\infty} \phi_N^2\left(y+\frac{x}{N}\right)\phi_N^2(y).$$

Recall that $\lim_{N\to\infty}\phi_N(x) = \phi(x)$ for almost every $x\in\mathbb{R}^d$. Therefore, it suffices to show that

$$(4.16) \qquad \limsup_{N\to\infty}\left|\phi_N^2(y) - \phi_N^2\left(y+\frac{x}{N}\right)\right| = 0 \qquad \text{for every } x,y\in\mathbb{R}^d.$$

Fix $x,y\in\mathbb{R}^d$. For a while, we write $z$ instead of $\frac{x}{N}$. Estimate

$$(4.17) \quad \begin{aligned} |\phi_N^2(y) - \phi_N^2(y+z)| &\leq \frac{1}{N}\sum_{i=1}^N |h_i^{(N)}(y) - h_i^{(N)}(y+z)| 2\|h_i\|_\infty \\ &\leq \frac{2C_d}{N}\sum_{i=1}^N (\lambda_i^{(N)})^{d/4}|h_i^{(N)}(y) - h_i^{(N)}(y+z)|. \end{aligned}$$

In the following we denote by $C$ a generic positive constant, which may change its value from appearance to appearance and does not depend on $i$ nor on $N$. Fix $i\in\{1,\dots,N\}$. By $U_i^{(N)}$, we denote the potential in (4.4) with $v$ replaced by $v_N$, that is,

$$U_i^{(N)}(w) = W(w) + \sum_{j\neq i} V_N h_j^2(w) - \lambda_i^{(N)},$$

where $V_N$ is the operator $V$ with $v$ replaced by $v_N$. Note that

$$(4.18) \quad \begin{aligned} \left\|\sum_{j\neq i} V_N h_j^2\right\|_\infty &\leq \frac{1}{N}\sum_{j=1}^N \sup_{x\in\mathbb{R}^d} \int_{\mathbb{R}^d} dy\, N^d v(N|x-y|) h_j^2(y) \\ &\leq \frac{C_d^2}{N}\sum_{j=1}^N (\lambda_j^{(N)})^{d/2}\int_{\mathbb{R}^d} v(|z|)\, dz \leq CN^{(d-2)/2}. \end{aligned}$$

Recall (4.2) and estimate

$$(4.19) \quad \begin{aligned} \lambda_i^{(N)} &\leq \inf_{\widetilde{h}_i\in H^1(\mathbb{R}^d),\|\widetilde{h}_i\|_2=1}[\|\nabla\widetilde{h}_i\|_2^2 + \langle W,\widetilde{h}_i^2\rangle + CN^{(d-2)/2}] \\ &\leq CN^{(d-2)/2}. \end{aligned}$$

We write $T := T_1(y)$ for the exit time of a Brownian motion $(B_t)_{t \geq 0}$ from the ball $B := B_1(y)$. By [3], Theorem A.4.1, Lemma A.4.5, we have

$$
\begin{aligned}
(4.20) \qquad h_i^{(N)}(w) &= \mathbb{E}_w[h_i^{(N)}(B_T)] \\
&\quad - \mathbb{E}_w\left[\int_0^T U_i^{(N)}(B_s) h_i^{(N)}(B_s)\,ds\right], \qquad w \in B.
\end{aligned}
$$

The first term is harmonic. By the mean-value property, we see that, denoting by $S_{y,z}$ the symmetric difference $B_{1/4}(y) \triangle B_{1/4}(y+z)$,

$$
\begin{aligned}
(4.21) \qquad |\mathbb{E}_y[h_i^{(N)}(B_T)] - \mathbb{E}_{y+z}[h_i^{(N)}(B_T)]| &\leq \frac{1}{|B_{1/4}(y)|} \int_{S_{y,z}} |\mathbb{E}_w[h_i^{(N)}(B_T)]|\,dw \\
&\leq \frac{\|h_i^{(N)}\|_\infty}{|B_{1/4}(y)|} |S_{y,z}| \\
&\leq C(\lambda_i^{(N)})^{d/4}|z|.
\end{aligned}
$$

The second term on the right-hand side of (4.20) may be written as $\widetilde{h}_i^{(N)}(w) + \widehat{h}_i^{(N)}(w)$, where

$$
(4.22) \qquad \widetilde{h}_i^{(N)}(w) := \int_B G(w, \widetilde{w}) U_i^{(N)}(\widetilde{w}) h_i^{(N)}(\widetilde{w})\,d\widetilde{w},
$$

$$
(4.23) \qquad \widehat{h}_i^{(N)}(w) := \int_B (G_B - G)(w, \widetilde{w}) U_i^{(N)}(\widetilde{w}) h_i^{(N)}(\widetilde{w})\,d\widetilde{w},
$$

and $G$ and $G_B$ are the Green's function in $\mathbb{R}^d$ and in the ball $B$, respectively. (For explicit formulas for $G_B$, see [8], Section 2.3, e.g.) By [21], Theorem 10.2(ii), for $p = d$ and any $\alpha \in (0, 1)$, we have

$$
|\widetilde{h}_i^{(N)}(y) - \widetilde{h}_i^{(N)}(y - z)| \leq C|z|^\alpha \|U_i^{(N)} h_i^{(N)} \mathbb{1}_B\|_p |B|^{(2-\alpha)/d - 1/p}.
$$

By (4.18) and (4.19), $\sup_{w \in B} |U_i^{(N)}(w)| \leq CN^{(d-2)/2}$. Hence, we can further estimate

$$
(4.24) \qquad |\widetilde{h}_i^{(N)}(y) - \widetilde{h}_i^{(N)}(y - z)| \leq C|z|^\alpha N^{(d-2)/2}(\lambda_i^{(N)})^{d/4}.
$$

Now we turn to $\widehat{h}_i^{(N)}$, which is harmonic in $B$, since $(G_B - G)(\cdot, \widetilde{w})$ is for any $\widetilde{w} \in B$. By [8], Theorem 2.5, $(G - G_B)(w, \widetilde{w}) = \mathbb{E}_w[G(B_T, \widetilde{w})]$ for any $w, \widetilde{w} \in B$. Like in (4.21), we use the mean-value property, recall that $\sup_{w \in B} |U_i^{(N)}(w)| \leq CN^{(d-2)/2}$, use Harnack's inequality and obtain, for any

$\alpha \in (0, 1)$,

$$
\begin{aligned}
|\widehat{h}_i^{(N)}(y) &- \widehat{h}_i^{(N)}(y-z)| \\
&\leq CN^{(d-2)/2}(\lambda_i^{(N)})^{d/4} \frac{1}{|B_{1/4}(y)|} \int_{S_{y,z}} dw \int_B d\widetilde{w} |(G - G_B)(w, \widetilde{w})| \\
&\leq CN^{(d-2)/2}(\lambda_i^{(N)})^{d/4} \frac{1}{|B_{1/4}(y)|} \int_{S_{y,z}} dw \int_B d\widetilde{w} \, \mathbb{E}_0[G(B_T, \widetilde{w})] \\
&\leq CN^{(d-2)/2}(\lambda_i^{(N)})^{d/4} \frac{|S_{y,z}|}{|B_{1/4}(y)|} \int_{B_4(0)} G(0, \widetilde{w}) \, d\widetilde{w} \\
&\leq CN^{(d-2)/2}(\lambda_i^{(N)})^{d/4} |z|^\alpha.
\end{aligned}
$$
(4.25)

Substituting the estimates (4.21), (4.24) and (4.25) in (4.17), we obtain

$$
|\phi_N^2(y) - \phi_N^2(y+z)| \leq C \frac{1}{N} \sum_{i=1}^N (\lambda_i^{(N)})^{d/2}(|z| + N^{(d-2)/2}|z|^\alpha + N^{(d-2)/2}|z|^\alpha).
$$

Recall that $z = x/N$ and, furthermore, that $\frac{1}{N} \sum_{i=1}^N \lambda_i^{(N)} \leq C$ and $\lambda_i^{(N)} \leq CN^{(d-2)/2}$. Picking $\alpha$ sufficiently close to 1, (4.14) easily follows. This finishes the proof of $\lim_{N\to\infty} \chi_N^{(\otimes)} = \chi_{\widetilde{\alpha}(v)}^{(\mathrm{GP})}$.

The proof of the convergence of $\phi_N$ toward the Gross–Pitaevskii minimizer $\phi$ (not only along subsequences) is done in a standard way as follows. Let $f \colon \mathbb{R}^d \to \mathbb{R}$ be a measurable bounded function, and fix some $\varepsilon \in \mathbb{R}$ with $|\varepsilon|$ small. Consider the trap potential $W_\varepsilon = W + \varepsilon f$. It satisfies all the required assumptions that we posed on $W$, with the possible exception of the nonnegativity. However, an examination of our above proofs shows that they work also for $W_\varepsilon$ in place of $W$. Let $g_N(\varepsilon)$ denote the variational formula in (1.11) with $W$ replaced by $W_\varepsilon$, then we know from the preceding that $\lim_{N\to\infty} g_N(\varepsilon)$ is equal to the Gross–Pitaevskii formula in (1.40) with $W$ replaced by $W_\varepsilon$. Furthermore, it is easily seen that $g_N$ is concave in a neighborhood of zero for any $N$. Hence, we know that also the derivatives $g_N'$ converge the $\varepsilon$-derivative of the Gross–Pitaevskii formula with $W$ replaced by $W_\varepsilon$. Clearly, $g_N'(0) = \langle f, \phi_N^2 \rangle$, where $\phi_N$ is as in the above proof of $\liminf_{N\to\infty} \chi_N^{(\otimes)} \geq \chi_{\alpha(v)}^{(\mathrm{GP})}$. Furthermore, the analogous derivative for the Gross–Pitaevskii formula is equal to $\langle f, \phi^2 \rangle$, where $\phi = \phi_{\widetilde{\alpha}(v)}^{(\mathrm{GP})}$ is the minimizer in (1.40) for $\alpha = \widetilde{\alpha}(v)$. By convergence of the derivatives, we have $\lim_{N\to\infty} \langle f, \phi_N^2 \rangle = \langle f, \phi^2 \rangle$. Indeed, for $\varepsilon > 0$, we have

$$
\limsup_{N\to\infty} \langle f, \phi_N^2 \rangle \leq \limsup_{N\to\infty} \frac{1}{\varepsilon}(g_N(\varepsilon) - g_N(0)) = \frac{1}{\varepsilon}(g(\varepsilon) - g(0)),
$$

and the right-hand side converges to $\langle f, \phi^2 \rangle$ as $\varepsilon \downarrow 0$. In order to see the reversed inequality, replace $f$ by $-f$ and $\varepsilon$ by $-\varepsilon$. This shows the convergence of $\phi_N^2(x)$ toward $\phi(x)^2$ both in the weak $L^1$-sense and in the sense of probability measures.

**5. Large deviations for the canonical ensemble model: Proof of Theorem 1.5.** In this section we prove large deviations statement on the canonical ensemble model in Theorem 1.5. We follow the same strategy as in the proof of Theorem 1.7 in Section 2. Hence, the main step is the proof of the following.

PROPOSITION 5.1 (Convergence of the logarithmic moment generating function). *For any* $f \in \mathcal{C}_{\mathrm{b}}(\Omega)$,

$$(5.1) \qquad \lim_{\beta \to \infty} \frac{1}{\beta} \log \mathbb{E}[e^{-\beta\langle \mathfrak{W} + \mathfrak{v} - f, \mu_\beta \rangle}] = -N\chi_N(f),$$

*where*

$$(5.2) \qquad \chi_N(f) = \frac{1}{N} \inf_{h \in H^1(\Omega)\,:\,\|h\|_2=1} [\|\nabla h\|_2^2 + \langle \mathfrak{W} + \mathfrak{v} - f, h^2 \rangle].$$

PROOF. See Sections 5.1 and 5.2 for the upper and lower bound, respectively. □

From Proposition 5.1, one derives Theorem 1.5 in the same way as Theorem 1.7 is derived from Proposition 2.1 in Section 2.4; we omit the details.

The assertion of Proposition 5.1 is classical and well known if the potential $\mathfrak{W} + \mathfrak{v} - f$ were bounded and continuous. But $\mathfrak{W}$ and $\mathfrak{v}$ are unbounded to $\infty$, and this is the additional technicality we are facing. Nevertheless, we feel that this problem and its solution are well known, but we could not find a precise reference. Hence, we provide a proof.

We shall employ eigenvalue expansions on large compact sub boxes to derive (5.1). For the upper bound, we use periodic boundary condition in the box, and for the lower bound, we use zero boundary condition. Let us remark here that, if $W$ and $v$ would have been assumed continuous in $\{W < \infty\}$, respectively in $\{v < \infty\}$, these arguments for the upper bound could have been also replaced by applications of Varadhan's lemma and the large deviations principles in Lemma 1.15. However, due to the degeneracy of $v$ at zero, we did not find any way to derive the lower bound in (5.1) via large deviations arguments in the soft-core case.

5.1. *Proof of the upper bound in* (5.1). Recall that $\Lambda_R = [-R, R]^d$. We shall divide the path space into the part on which the Brownian motion spends more than $(1-\eta)\beta$ time units in $\Lambda_R^N$ up to time $\beta$ and the remaining part where this is not satisfied. We will replace $\mu_\beta$ by its $\Lambda_R^N$-periodized version $\mu_{\beta,R}$ and control the error. We also replace the two potentials $W$ and $v$ by their cut-off versions $W_M = W \wedge M$ and $v_M = v \wedge M$, where $M > 0$ is large. Finally, we let $R \to \infty$, $\eta \downarrow 0$ and $M \to \infty$.

In the following we sketch only the proof because all the details can be found in the corresponding proof for the Hartree model. For the canonical ensemble model, it turns out the error estimates are simpler, in particular, the interaction part, because here only *one* Brownian motion is involved, and hence, there is only *one* time scale. We begin with

$$(5.3) \quad \mathbb{E}(e^{-\beta\langle\mathfrak{W}+\mathfrak{v}-f,\mu_\beta\rangle}) \leq \mathbb{E}(e^{-\beta\langle\mathfrak{W},\mu_\beta\rangle}\mathbb{1}\{\mu_\beta(\Lambda_R^N) \leq 1-\eta\})e^{-\beta\inf v}e^{\beta\|f\|_\infty} \\ + \mathbb{E}(e^{-\beta\langle\mathfrak{W}_M+\mathfrak{v}_M-f,\mu_\beta\rangle}\mathbb{1}\{\mu_\beta(\Lambda_R^N) > 1-\eta\}),$$

where $\mathfrak{W}_M$ and $\mathfrak{v}_M$ are defined as $\mathfrak{W}$ and $\mathfrak{v}$ in (1.3) with $W$ and $v$ replaced by $W \wedge M$ and $v \wedge M$, respectively. The first term on the right-hand side of (5.3) is easily further estimated, using that $W \geq 0$. Indeed, it is not bigger than $\exp\{-\beta(\eta\inf_{\Lambda_R^c} W + \inf v - \|f\|_\infty)\}$. For estimating the second term, we replace $\mu_\beta$ by its periodized version $\mu_{\beta,R}$ as in (1.45) with the box $\Lambda$ replaced by $\Lambda_R^N$. Now we can estimate the replacement error for the trap part as follows:

$$(5.4) \quad \begin{aligned} &\langle\mathfrak{W}_M,\mu_\beta\rangle - \langle\mathfrak{W}_M,\mu_{\beta,R}\rangle \\ &= \sum_{k\in\mathbb{Z}^{dN}\setminus\{0\}}\int_{\Lambda_R^N}\sum_{i=1}^N(W(x_i+2Rk_i)-W(x_i))\mu_\beta(d(x+2Rk)) \\ &\geq -M\mu_\beta((\Lambda_R^N)^c) \\ &\geq -\eta M. \end{aligned}$$

The replacement errors of the other parts of the potential are estimated in the same way:

$$\begin{aligned} &\langle\mathfrak{v}_M-f,\mu_\beta\rangle - \langle\mathfrak{v}_M-f,\mu_{\beta,R}\rangle \\ &\geq (-M+\inf v - 2\|f\|_\infty)\mu_\beta((\Lambda_R^N)^c) \\ &= -\eta(M-\inf v + 2\|f\|_\infty). \end{aligned}$$

Summarizing, we obtain

$$(5.5) \quad \begin{aligned} \text{l.h.s. of (5.1)} &\leq e^{-\beta(\eta\inf_{\Lambda_R^c} W + \inf v - \|f\|_\infty)} \\ &+ e^{2\beta\eta(M-\inf v + \|f\|_\infty)}\mathbb{E}(e^{-\beta\langle\mathfrak{W}_M+\mathfrak{v}_M-f,\mu_{\beta,R}\rangle}). \end{aligned}$$

Now we use an eigenvalue expansion to derive that

$$(5.6) \quad \limsup_{\beta\to\infty}\frac{1}{\beta}\log\mathbb{E}(e^{-\beta\langle\mathfrak{W}_M+\mathfrak{v}_M-f,\mu_{\beta,R}\rangle}) \leq -N\chi_N(R,M,f),$$

where

$$(5.7) \quad \begin{aligned} \chi_N(R,M,f) = \frac{1}{N}\inf\Big\{&\frac{1}{2}\int_{\Lambda_R^N\cap\Omega}|\nabla_R h(x)|^2\,dx + \langle\mathfrak{W}_M+\mathfrak{v}_M-f,h^2\rangle : \\ &h \in \mathcal{C}^\infty(\Lambda_R^N\cap\Omega), \|h\|_2 = 1\Big\}. \end{aligned}$$

Here $\nabla_R$ denotes the periodized gradient on $\Lambda_R^N$, that is, the one with periodic boundary condition. To derive (5.6), we let $(\lambda_k)_{k\in\mathbb{N}}$ and $(e_k)_{k\in\mathbb{N}}$ be the sequence of eigenvalues and an orthonormal basis of corresponding eigenfunctions of the operator $\Delta - \mathfrak{W}_M - \mathfrak{v}_M + f$ in $\Lambda_R^N \cap \Omega$ with periodic boundary condition. We may assume that $\lambda_1 = -N\chi_N(R, M, f)$ is the principal eigenvalue. Let $p_1^{(R)}(x, y)$ denote the transition probability function of the Brownian motion on the torus $\Lambda_R^N \cap \Omega$ and note that $p_1^{(R)}(x, y) \le C_R$ for any $x, y \in \Lambda_R^N$, for some $C_R < \infty$. For simplicity, we assume that $\nu(\Lambda_R^N) = 1$. Then we can estimate

$$\mathbb{E}(e^{-\beta\langle \mathfrak{W}_M + \mathfrak{v}_M - f, \mu_{\beta,R}\rangle})$$

$$\le e^{-\inf v + \|f\|_\infty} \int \nu(dx)\mathbb{E}_x(e^{-\int_1^\beta [\mathfrak{W}_M + \mathfrak{v}_M - f](B_s)\,ds})$$

$$= e^{-\inf v + \|f\|_\infty} \int_{\Lambda_R^N} dy \int \nu(dx) p_1^{(R)}(x, y) \mathbb{E}_y(e^{-(\beta-1)\langle \mathfrak{W}_M + \mathfrak{v}_M - f, \mu_{\beta-1,R}\rangle})$$

$$\le C \sum_{k\in\mathbb{N}} e^{(\beta-1)\lambda_k} \langle e_k, \mathbb{1}\rangle^2,$$

where $C > 0$ does not depend on $\beta$, and $\langle \cdot, \cdot \rangle$ is the inner product on $L^2(\Lambda_R^N)$. Now use Parseval's identity to continue

$$\mathbb{E}(e^{-\beta\langle \mathfrak{W}_M + \mathfrak{v}_M - f, \mu_{\beta,R}\rangle}) \le e^{(\beta-1)\lambda_1} \sum_{k\in\mathbb{N}} \langle e_k, \mathbb{1}\rangle^2$$

$$= e^{(\beta-1)\lambda_1} \|\mathbb{1}\|^2,$$

where $\|\cdot\|$ denotes the $L^2$ norm in $\Lambda_R^N$. From this, (5.6) directly follows. From (5.5) and (5.6), we obtain

(5.8)
$$\begin{aligned}
\text{l.h.s. of } &(5.1) \\
&\le -\min\Big\{ \eta \inf_{\Lambda_R^c} W + \inf v - \|f\|_\infty, \\
&\qquad\qquad -2\eta[M - \inf v + \|f\|_\infty] + N\chi_N(R, M, f) \Big\}.
\end{aligned}$$

Now we let $R \to \infty$, then $\eta \downarrow 0$ and finally $M \to \infty$ on the right-hand side of (5.8). Note that

(5.9)
$$\liminf_{M\to\infty} \liminf_{R\to\infty} \chi_N(R, M, f) \ge \chi_N(f).$$

The proof of (5.9) is standard, we do not carry it out here; note that the proof of (2.14) is similar and even technically more difficult.

Because of our assumption on $W$ in (1.1), $\lim_{R\to\infty} \inf_{\Lambda_R^c} W = \infty$. Hence, we obtain that the left-hand side of (5.1) is not bigger than $N\chi_N(f)$. This completes the proof of the upper bound in (5.1).

5.2. *Proof of the lower bound in* (5.1). Now we prove the lower bound in (5.1). Like in the proof of the upper bound, we rely on an eigenvalue expansion. We recall $\Lambda_R = [-R, R]^d$ and the set $U_\eta$ from (1.5). Fix some $\eta > a$ [recall (1.2)]. We estimate the expectation on the left-hand side of (5.1) by imposing zero boundary condition in a certain compact set:

$$\mathbb{E}[e^{-\beta \langle \mathfrak{W} + \mathfrak{v} - f, \mu_\beta \rangle}] \geq \mathbb{E}[e^{-\beta \langle \mathfrak{W} + \mathfrak{v} - f, \mu_\beta \rangle} \mathbb{1}\{\operatorname{supp}(\mu_\beta) \subset \Lambda_R^N \cap \widetilde{U}_\eta\}],$$

where we put $\widetilde{U}_\eta = U_\eta$ in the hard-core case and $\widetilde{U}_\eta = \mathbb{R}^{dN}$ in the soft-core case. Hence, in the hard-core case, we may replace $\mathfrak{W}$ and $\mathfrak{v}$ by $\mathfrak{W}_M$ and $\mathfrak{v}_M$ [for this notation, see below (5.3)], respectively, where $M > 0$ depends only on $\eta$ and $R$. In the soft-core case, the potential $\mathfrak{W} + \mathfrak{v} - f$ is locally integrable in the box $\Lambda_R$ and bounded from below.

Consider the linear operator $\Delta - \mathfrak{W} - \mathfrak{v}_f$ on the space $L^2(\Lambda_R^N \cap \widetilde{U}_\eta)$ with zero boundary condition. Standard results imply that this operator possesses a compact resolvent. Hence, a Fourier expansion in terms of the eigenvalues of this operator yields that

(5.10)
$$\lim_{\beta \to \infty} \frac{1}{\beta} \log \mathbb{E}[e^{-\beta \langle \mathfrak{W} + \mathfrak{v} - f, \mu_\beta \rangle} \mathbb{1}\{\operatorname{supp}(\mu_\beta) \subset \Lambda_R^N \cap \widetilde{U}_\eta\}]$$
$$= -N \chi_N(R, \eta, f),$$

where

(5.11)
$$\chi_N(R, \eta, f)$$
$$= \frac{1}{N} \inf_{h \in H^1(\mathbb{R}^{dN}), \operatorname{supp}(h) \subset \Lambda_R^N \cap \widetilde{U}_\eta, \|h\|_2 = 1} \left( \frac{1}{2} \|\nabla h\|_2^2 + \langle \mathfrak{W} + \mathfrak{v} - f, h^2 \rangle \right)$$

is the principal eigenvalue of the above operator. The proof of (5.10) is similar to the one of (2.31) and technically much easier; we omit the details.

Now we let $R \to \infty$ and $\eta \downarrow a$ and see easily that

$$\lim_{R \to \infty} \lim_{\eta \downarrow a} \chi_N(R, \eta, f) = \chi_N(f).$$

This completes the proof of the lower bound in (5.1).

**Acknowledgments.** Wolfgang König wishes to thank the DFG for awarding a Heisenberg grant. All authors gratefully acknowledge the hospitality of the Dublin Institute for Advanced Studies where a substantial part of this work was carried out.

S. ADAMS
MAX-PLANCK INSTITUTE FOR
   MATHEMATICS IN THE SCIENCES
INSELSTRASSE 22-26
D-04103 LEIPZIG
GERMANY
AND
SCHOOL OF THEORETICAL PHYSICS
DUBLIN INSTITUTE FOR ADVANCED STUDIES
10 BURLINGTON ROAD
DUBLIN 4
IRELAND
E-MAIL: adams@mis.mpg.de
URL: www.mis.mpg.de/sm/homepages/adams.html

J.-B. BRU
FACHBEREICH MATHEMATIK UND INFORMATIK
JOHANNES-GUTENBERG-UNIVERSITÄT MAINZ
STAUDINGERWEG 9
D-55099 MAINZ
GERMANY
E-MAIL: jbbru@mathematik.uni-mainz.de
URL: www.mathematik.uni-mainz.de/Members/jbbru

W. KÖNIG
MATHEMATISCHES INSTITUT
UNIVERSITÄT LEIPZIG
AUGUSTUSPLATZ 10/11
D-04109 LEIPZIG
GERMANY
E-MAIL: koenig@math.uni-leipzig.de
URL: www.math.uni-leipzig.de/~koenig/